\documentclass[reqno]{amsart}
\usepackage{amsmath,amsfonts,amssymb,xcolor,enumitem,upgreek}
\usepackage[english]{babel}
\usepackage[colorlinks=true,linkcolor=red!60!black,citecolor=green!60!black]{hyperref}

\newtheorem{theorem}{Theorem}[section]
\newtheorem{lemma}[theorem]{Lemma}
\newtheorem{proposition}[theorem]{Proposition}
\newtheorem{corollary}[theorem]{Corollary}

\theoremstyle{definition}

\newcommand{\w}{\omega}

\numberwithin{equation}{section}

\begin{document}
\title[Generalized Pitchfork Bifurcation]{Generalized Pitchfork Bifurcations
\\in D-Concave Nonautonomous Scalar
\\Ordinary Differential Equations}
\author[J. Due\~{n}as]{Jes\'{u}s Due\~{n}as}
\author[C. N\'{u}\~{n}ez]{Carmen N\'{u}\~{n}ez}
\author[R. Obaya]{Rafael Obaya}
\address{Departamento de Matem\'{a}tica Aplicada, Universidad de Va\-lladolid, Paseo Prado de la Magdalena 3-5, 47011 Valladolid, Spain. The authors are members of IMUVA: Instituto de Investigaci\'{o}n en Matem\'{a}ticas, Universidad de Valladolid.}
\email[J.~Due\~{n}as]{jesus.duenas@uva.es}
\email[C.~N\'{u}\~{n}ez]{carmen.nunez@uva.es}
\email[R.~Obaya]{rafael.obaya@uva.es}
\thanks{All the authors were supported by Ministerio de Ciencia, Innovaci\'{o}n y Universidades (Spain) under project RTI2018-096523-B-I00 and by Universidad de Valladolid under project PIP-TCESC-2020. J. Due\~{n}as was supported by Ministerio de Universidades (Spain) under programme FPU20/01627.}
\date{}
\begin{abstract}
The global bifurcation diagrams for two different one-parametric perturbations ($+\lambda x$ and $+\lambda x^2$) of a dissipative scalar nonautonomous ordinary differential equation $x'=f(t,x)$ are described assuming that 0 is a constant solution, that $f$ is recurrent in $t$, and that its first derivative with respect to $x$ is a strictly concave function. The use of the skewproduct formalism allows us to identify  bifurcations with changes in the number of minimal sets and in the shape of the global attractor. In the case of perturbation $+\lambda x$, a so-called {\em global generalized pitchfork bifurcation} may arise, with the particularity of lack of an analogue in autonomous dynamics. This new bifurcation pattern is extensively investigated in this work.
\end{abstract}
\renewcommand{\subjclassname}{\textup{2020} Mathematics Subject Classification}

\maketitle
\section{Introduction}\label{sec:1introduction}
The interest that the description of nonautonomous bifurcation patterns
arouses in the scientific community has increased significantly in recent years,
as evidenced by the works \cite{alob3}, \cite{anagjager1}, \cite{anjk}, \cite{brbh}, \cite{dno1},
\cite{fabri1}, \cite{joma}, \cite{johnson1}, \cite{kloeden1},
\cite{langa1}, \cite{nunezobaya,nuob10}, \cite{potz2,potz4}, \cite{rasmussen2,rasmussen1}, \cite{refj}, and references therein. This paper constitutes an extension of the work initiated
in \cite{dno1}, were we describe several possibilities for the global bifurcation
diagrams of certain types of one-parametric variations of a coercive equation. We
make use of the skewproduct formalism, which allows us to understand bifurcations
as changes in the number of minimal sets and in the shape of the global attractor,
which of course give rise to substantial changes in the global dynamics.

Let us briefly describe the skewproduct formalism. Standard boundedness
and regularity conditions ensure that the hull $\Omega$ of a continuous map
$f_0\colon\mathbb{R}\times\mathbb{R}\rightarrow\mathbb{R}$, defined as the closure
of the set of time-shifts $\{f_0{\cdot}t\colon\, t \in \mathbb{R}\}$ in a suitable
topology of $C(\mathbb{R}\times\mathbb{R},\mathbb{R})$, is a compact metric space,
and that the map
$\sigma\colon\mathbb{R}\times\Omega\rightarrow\Omega,\;(t,\omega)\mapsto
\omega{\cdot}t$ (where, as in the case of $f_0$,
$(\omega{\cdot}t)(s,x)=\omega(t+s,x)$) defines a global continuous flow. The
continuous function $f(\omega,x)=\omega(0,x)$ provides the family of equations
\begin{equation}\label{1.primera}
 x'=f(\omega{\cdot}t,x)\,,\qquad \omega\in\Omega\,,
\end{equation}
which includes $x'=f_0(t,x)$: it corresponds to $\omega_0=f_0\in\Omega$.
When, in addition, $f_0$ satisfies some properties of recurrence in time,
the flow $(\Omega,\sigma)$ is minimal, which means that $\Omega$ is
the hull of any of its elements. If $u(t,\omega,x)$ denotes the solution
of \eqref{1.primera}$_\omega$ with $u(0,\omega,x)=x$, then
$\tau\colon\mathbb{R}\times\Omega\times\mathbb{R}\rightarrow\Omega\times\mathbb{R}$,
$(t,\omega,x)\mapsto(\omega{\cdot}t,u(t,\omega,x))$ defines a local
flow on $\Omega\times\mathbb{R}$ of skewproduct type: it projects over
the flow $(\Omega, \sigma)$. If $f_0$ is coercive with respect to
$x$ uniformly in $t\in\mathbb R$, so is $f$ uniformly in $\omega\in\Omega$,
and this ensures the existence of the global attractor and of at least one
minimal compact $\tau$-invariant subset of $\Omega\times\mathbb R$. In the simplest nonautonomous cases, the minimal subsets are
(hyperbolic or nonhyperbolic) graphs of continuous functions, and thus they play
the role performed by the critical points of an autonomous equation;
but there are cases in which both the shape of a minimal set and the
dynamics on it are extremely complex, without autonomous analogues,
and therefore impossible bifurcation scenarios for a
autonomous equation can appear in the nonautonomous setting.

So, we take as starting points a (global) continuous minimal flow $(\Omega,\sigma)$
and a continuous map $f\colon\Omega\times\mathbb R\to\mathbb R$, assume that $f$ is
coercive in $x$ uniformly on $\Omega$, and define the dissipative flow $\tau$.
Throughout this paper, we also assume
that the derivatives $f_x$ and $f_{xx}$ globally exist and are jointly continuous on
$\Omega\times\mathbb R$, as well as the fundamental property of strict concavity of $f_x$
with respect to $x$: {\em d-concavity}.
Not all these conditions are in force to obtain the results of
\cite{dno1}, but, for simplicity, we also assume them all to summarize part of the properties
there proved.

The first goal in \cite{dno1} is to describe the possibilities for the global
$\mu$-bifurcation diagram of the one-parametric family $x'=f(\omega{\cdot}t,x)+\mu$,
with global attractor $\mathcal A_\mu$. In particular, it is proved that, if there exist
three minimal sets for a value $\mu_0\in\mathbb{R}$ of the parameter,
then: $\mathcal{A}_\mu$ contains three (hyperbolic) minimal sets if and only if
$\mu$ belongs to a nondegenerate interval $(\mu_-,\mu_+)$;
the two upper (resp.~lower) minimal sets collide on a residual invariant subset of
$\Omega$ when $\mu\downarrow\mu_-$ (resp.~$\mu\uparrow\mu_+$);
and $\mathcal{A}_\mu$ reduces to the (hyperbolic) graph of a continuous map on $\Omega$ if
$\mu\notin[\mu_-,\mu_+]$. That is, the global bifurcation diagram presents at $\mu_-$
and $\mu_+$ two {\em local saddle-node bifurcation points}
of minimal sets and two points of discontinuity
of $\mathcal{A}_{\mu}$: is the nonautonomous analogue of the bifurcation diagram of $x'=-x^3+x+\mu$.

A second type of perturbation is considered in \cite{dno1}, namely $x'=f(\omega{\cdot}t,x)+\lambda x$,
with global attractor $\mathcal A_\lambda$, under the additional assumption $f({\cdot},0)\equiv 0$.
Now, $\mathcal M_0=\Omega\times\{0\}$ is a minimal set for all $\lambda$, and its hyperbolicity
properties are determined by the Sacker and Sell spectrum $[-\lambda_+,-\lambda_-]$
of the map $\omega\mapsto f_x(\omega,0)$.
Two possible global bifurcation diagrams are described, and some conditions ensuring their occurrence
are given. The first one is the {\em classical global pitchfork bifurcation diagram}, with unique
bifurcation point $\lambda_+$: $\mathcal M_0$ is the unique minimal set for $\lambda\le\lambda_+$,
and two more (hyperbolic) minimal sets occur for $\lambda>\lambda_+$, which collide with $\mathcal M_0$ as
$\lambda\downarrow\lambda_+$. An autonomous analogue is the diagram of $x'=-x^3+\lambda$.
The second one is the {\em local saddle-node and transcritical bifurcation diagram}, with a local
saddle node bifurcation of minimal sets at a point $\lambda_0<\lambda_-$ and a so-called {\em generalized} transcritical bifurcation of minimal sets around $\mathcal M_0$.
We will describe this diagram in detail in the next pages, pointing out now the most remarkable fact:
$\mathcal M_0$ collides with another (hyperbolic) minimal set as $\lambda\uparrow\lambda_-$ and
as $\lambda\downarrow\lambda_+$, and it is the unique minimal set lying on a band $\Omega\times[-\rho,\rho]$
for a $\rho>0$ if $\lambda\in[\lambda_-,\lambda_+]$. This local transcritical bifurcation becomes
classical if $\lambda_-=\lambda_+$, being $x'=-x^3+2x^2+\lambda x$ an autonomous example of this situation.

This analysis of the family $x'=f(\omega{\cdot}t,x)+\lambda x$ initiated in \cite{dno1} is
far away to be complete. The goal of this paper is to describe {\em all\/} the possibilities for its global
bifurcation diagram. Besides the two described ones, only a third situation may arise:
a {\em global generalized pitchfork bifurcation diagram}, just possible when $\lambda_-<\lambda_+$.
It is characterized by the existence of two bifurcation points, $\lambda_0\in[\lambda_-,\lambda_+)$ and
$\lambda_+$: $\mathcal M_0$ is the unique minimal set for $\lambda<\lambda_0$, there are two of them for
$\lambda\in(\lambda_0,\lambda_+]$, and there are three for $\lambda>\lambda_+$. The lack of an autonomous analogue
raises a nontrivial question: does this bifurcation diagram correspond to some actual family? We also
answer it, explaining how to construct nonautonomous patterns fitting at each one of the
described possibilities. Furthermore, we prove that, given $\lambda_-<\lambda_+$,
any $\lambda_0\le\lambda_+$ is the first bifurcation point of a suitable family
$x'=g_{\lambda_0}(\omega{\cdot}t,x)+\lambda x$ with Sacker and Sell spectrum of $(g_{\lambda_0})_x({\cdot},0)$
given by $[-\lambda_+,-\lambda_-]$, and that the three possible diagrams actually occur: they correspond
to $\lambda_0<\lambda_-$, $\lambda_0=\lambda_+$ and $\lambda_0\in[\lambda_-,\lambda_+)$. As a tool to
prove of this last result, we analyze the bifurcation possibilities for a new one-parametric family,
namely $x'=f(\omega{\cdot}t,x)+\xi x^2$. In order not to lengthen this introduction too much, we
omit here the (self-interesting) description of the bifurcation possibilities for this case, and refer
the reader to Section \ref{sec:7anotherbifurcation}.

These are the main results of this paper, which presents more detailed descriptions in some
particular cases. Its contents are organized in five sections. Section \ref{sec:2preliminiaries} contains
the basic notions and properties required to start with the analysis. Section \ref{sec:3bifurcationtheorem}
is devoted to the description of the three mentioned possibilities for the bifurcation diagrams of
$x'=f(\omega{\cdot}t,x)+\lambda x$. In Section \ref{sec:5criteriacubic}, we focus on the case of a
cubic polynomial $f(\w,x)=-a_3(\w)x^3+a_2(\w)x^2+a_1(\w)$ with strictly positive $a_3$, and show
how some suitable properties of the coefficients $a_1, a_2$ and $a_3$ and some factible relations
among them either preclude or ensure each one of the three different bifurcation diagrams. Section \ref{sec:6generalframework}
extends these results to more general functions $f(\w,x)=(-a_3(\w)+h(\w,x))x^3+a_2(\w)x^2+a_1(\w)$, describing
in this way other patterns fitting each one of the possibilities. And Section \ref{sec:7anotherbifurcation}
begins with the description of the casuistic for the bifurcation diagrams of $x'=f(\omega{\cdot}t,x)+\xi x^2$
to conclude with the consequence mentioned at the end of the previous paragraph.

\section{Preliminaries}\label{sec:2preliminiaries}
Throughout the paper, the map $\sigma\colon\mathbb{R}\times\Omega\rightarrow\Omega$, $(t,\omega)\mapsto \sigma_t(\omega)=\omega{\cdot}t$ defines a global continuous flow on a compact metric space $\Omega$, and we assume that the flow $(\Omega, \sigma)$ is minimal, that is, that every $\sigma$-orbit is dense in $\Omega$. This paper will be focused on describing the bifurcation diagrams of simple parametric variations of the family
\begin{equation}\label{eq:generalequation}
x'=f(\omega{\cdot}t,x)\,,\quad\omega\in\Omega\,,
\end{equation}
where $f\colon\Omega\times\mathbb{R}\rightarrow\mathbb{R}$ is assumed to be jointly continuous, $f_x$ and $f_{xx}$ are supposed to exist and to be jointly continuous (which we represent as $f\in C^{0,2}(\Omega\times\mathbb{R},\mathbb{R})$), and $f(\omega,0)=0$ for all $\omega\in\Omega$ (that is, $x\equiv0$ solves the equation). If only $f$ and $f_x$ are assumed to exist and to be jointly continuous, then we shall say that $f\in C^{0,1}(\Omega\times\mathbb{R},\mathbb{R})$. Additional coercivity and concavity properties will be assumed throughout the paper. In Section~\ref{sec:5criteriacubic}, we focus on the case in which $f(\omega,x)$ is a cubic polynomial in the state variable $x$ with strictly negative cubic coefficient.

We develop our bifurcation theory through the skewproduct formalism: as explained in the Introduction, our bifurcation analysis studies the variations on the global attractors and on the number and structure of minimal sets for the corresponding parametric family of skewproduct flows. In the next subsections, we summarize the most basic concepts and some basic results required in the formulations and proofs of our results. The interested reader can find in Section~2 of \cite{dno1} more details on these matters, as well as a suitable list of references.
\subsection{Scalar skewproduct flow}\label{subsec:skewproductflow}
We define the \emph{local skewproduct flow} on $\Omega\times\mathbb{R}$ induced by \eqref{eq:generalequation} as
\begin{equation}\label{eq:flujotau}
\tau\colon\mathcal{U}\subseteq \mathbb{R}\times\Omega\times\mathbb{R}\rightarrow\Omega\times\mathbb{R}\,,\quad (t,\omega,x_0)\mapsto (\omega{\cdot}t,u(t,\omega,x_0))\,,
\end{equation}
where $\mathcal{I}_{\omega,x_0}\rightarrow\mathbb{R}$, $t\mapsto u(t,\omega,x_0)$ is the maximal solution of $\eqref{eq:generalequation}$ with initial value $u(0,\omega,x_0)=x_0$, and $\mathcal{U}=\bigcup_{(\omega,x_0)\in\Omega\times\mathbb{R}}(\mathcal{I}_{\omega,x_0}\times\{(\omega,x_0)\})$. The flow $\tau$ projects on $(\Omega,\sigma)$, which is called its \emph{base flow}. It is known that $u$ satisfies the \emph{cocycle property} $u(t+s,\omega,x_0)=u(t,\omega{\cdot}s,u(s,\omega,x_0))$ whenever the right-hand term is defined; and, clearly, the flow $\tau$ is fiber-monotone: if $x_1<x_2$, then $u(t,\omega,x_1)<u(t,\omega,x_2)$ for any $t$ in the common interval of definition of both solutions.

\subsection{Functions of bounded primitive}\label{subsec:2boundedprimitive} Throughout this paper, the space of continuous functions from $\Omega$ to $\mathbb{R}$ will be represented by $C(\Omega)$, the subspace of functions $a\in C(\Omega)$ such that $\int_\Omega a(\omega)\, dm=0$ for all $m\in\mathfrak{M}_\mathrm{erg}(\Omega,\sigma)$ will be represented by $C_0(\Omega)$, the subspace of functions $a\in C(\Omega)$ such that the map $t\mapsto a_\omega(t)=a(\omega{\cdot}t)$ is continuously differentiable on $\mathbb{R}$ will be represented by $C^1(\Omega)$ (in this case we shall represent $a'(\omega)=a_\omega'(0)$), and the subspace of functions $a\in C(\Omega)$ with continuous primitive, that is, such that there exists $b\in C^1(\Omega)$ with $b'=a$, will be represented by $CP(\Omega)$.

It is frequent to refer to a function $a\in CP(\Omega)$ as ``with bounded primitive''. Let us explain briefly the reason. Recall that $(\Omega,\sigma)$ is minimal. Then, $a\in CP(\Omega)$ if and only if there exists $\omega_0\in\Omega$ such that the map $a_0\colon\mathbb{R}\rightarrow\mathbb{R}$, $t\mapsto a(\omega_0{\cdot}t)$ has a bounded primitive $b_0(t)=\int_0^t a_0(s)\, ds$, in which case this happens for all $\omega\in\Omega$ (see e.g. Lemma~2.7 of \cite{johnson3} or Proposition~A.1 of \cite{jnot}).


It follows from Birkoff's Ergodic Theorem that $CP(\Omega)\subseteq C_0(\Omega)$. It is well known that $CP(\Omega)= C_0(\Omega)$ if $(\Omega,\sigma)$ is periodic, and Lemma~5.1 of \cite{cot1} ensures that $CP(\Omega)$ is a dense subset of $C_0(\Omega)$ of first category in $C_0(\Omega)$ if $(\Omega,\sigma)$ is a minimal aperiodic flow. In addition, $C^1(\Omega)$ is dense in $C(\Omega)$ (see Section~2 of \cite{sol1}).

\subsection{Equilibria and semiequilibria}\label{subsec:2equilibria} Let $\beta:\Omega\rightarrow\mathbb{R}$ be such that $u(t,\omega,\beta(\omega))$ is defined for all $\omega\in\Omega$ and $t\geq0$ (resp. $t\leq0$). We shall say that $\beta$ is a $\tau$-\emph{equilibrium} if $\beta(\omega{\cdot}t)=u(t,\omega,\beta(\omega))$ for all $\omega\in\Omega$ and $t\in\mathbb{R}$, a $\tau$-\emph{subequilibrium} (resp. \emph{time-reversed $\tau$-subequilibrium}) if $\beta(\omega{\cdot}t)\leq u(t,\omega,\beta(\omega))$ for all $\omega\in\Omega$ and $t\geq 0$ (resp. $t\leq0$), and a $\tau$-\emph{superequilibrium} (resp. \emph{time-reversed $\tau$-superequilibrium}) if $\beta(\omega{\cdot}t)\geq u(t,\omega,\beta(\omega))$ for all $\omega\in\Omega$ and all $t\geq 0$ (resp. $t\leq0$). Superequilibria and subequilibria are generally denominated \emph{semiequilibria}. We will frequently omit the reference to $\tau$ if there is no risk of confusion.

Given a Borel measure $m$ on $\Omega$, we shall say that $\beta\colon\Omega\rightarrow\mathbb{R}$ is $m$-\emph{measurable} if it is measurable with respect to the $m$-completion of the Borel $\sigma$-algebra, and we shall say that $\beta\colon\Omega\rightarrow\mathbb{R}$ is $C^1$ \emph{along the base orbits} if $t\mapsto\beta_\omega(t)=\beta(\omega{\cdot}t)$ is $C^1$ on $\mathbb{R}$ for all $\omega\in\Omega$, in which case we represent $\beta'(\omega)=\beta'_\omega(0)$. Note that any $\tau$-equilibrium is $C^1$ along the base orbits. We shall say that $\beta\colon\Omega\rightarrow\mathbb{R}$ is a \emph{semicontinuous equilibrium} (resp.~\emph{semiequilibrium}) if it is an equilibrium (resp.~semiequilibrium) and a bounded semicontinuous map. A \emph{copy of the base for the flow $\tau$} is the graph of a continuous $\tau$-equilibrium.

We shall say that a $\tau$-superequilibrium (resp. time-reversed $\tau$-superequilibrium) $\beta\colon\Omega\rightarrow\mathbb{R}$ is \emph{strong} if there exists $s_*>0$ (resp. $s_*<0$) such that $\beta(\omega{\cdot}s_*)>u(s_*,\omega,\beta(\omega))$ for all $\omega\in\Omega$, and we shall say that a $\tau$-subequilibrium (resp. time-reversed $\tau$-subequilibrium) is \emph{strong} if there exists $s_*>0$ (resp. $s_*<0$) such that $\beta(\omega{\cdot}s_*)<u(s_*,\omega,\beta(\omega))$ for all $\omega\in\Omega$. Since $(\Omega,\sigma)$ is minimal, Proposition~4.7 of \cite{nno1} and its analogue in the time-reversed cases show that, if $\beta$ is a semicontinuous strong $\tau$-superequilibrium (resp. time-reversed $\tau$-superequilibrium), then there exist $e_0>0$ and $s^*>0$ (resp. $s^*<0$) such that $\beta(\omega)\geq u(s^*,\omega{\cdot}(-s^*),\beta(\omega{\cdot}(-s^*)))+e_0$ for all $\omega\in\Omega$, and if $\beta$ is a semicontinuous strong $\tau$-subequilibrium (resp. time-reversed $\tau$-subequilibrium), then there exist $e_0>0$ and $s^*>0$ (resp. $s^*<0$) such that $\beta(\omega)\leq u(s^*,\omega{\cdot}(-s^*),\beta(\omega{\cdot}(-s^*)))-e_0$ for all $\omega\in\Omega$.

Let $\beta\colon\Omega\rightarrow\mathbb{R}$ be $C^1$ along the base orbits. The map $\beta$ shall be said to be a \emph{global upper} (resp.~\emph{lower}) \emph{solution} of \eqref{eq:generalequation} if $\beta'(\omega)\geq f(\omega,\beta(\omega))$ (resp.~$\beta'(\omega)\leq f(\omega,\beta(\omega))$) for every $\omega\in\Omega$, and to be \emph{strict} if the previous inequalities are strict for all $\omega\in\Omega$. Some comparison arguments prove the following facts (see Sections 3 and 4 of \cite{nno1}): if every forward $\tau$-semiorbit is globally defined, then $\beta$ is a $\tau$-superequilibrium (resp.~$\tau$-subequilibrium) if and only if is a global upper (resp.~lower) solution of \eqref{eq:generalequation}, and it is strong as superequilibrium (resp. subequilibrium) if it is strict as global upper (resp. lower) solution. Analogously, if every backward $\tau$-semiorbit is globally defined, then $\beta$ is a time-reversed $\tau$-subequilibrium (resp.~time-reversed $\tau$-superequilibrium) if and only if it is a global upper (resp.~lower) solution of \eqref{eq:generalequation}, and it is strong as time-reversed subequilibrium (resp. time-reversed superequilibrium) if it is strict as global upper (resp. lower) solution.
\subsection{Minimal sets, coercivity and global attractor}\label{subsec:2compactcoerciveattractor}
A set $\mathcal{K}\subseteq\Omega\times\mathbb{R}$ is $\tau$-\emph{invariant} if it is composed by globally defined $\tau$-orbits, and it is \emph{minimal} if it is compact, $\tau$-invariant and it does not contain properly any other compact $\tau$-invariant set. Let us recall some properties of compact $\tau$-invariant sets and minimal sets for the local skewproduct flow $(\Omega\times\mathbb{R},\tau)$ over a minimal base $(\Omega,\sigma)$. Let $\mathcal K\subset\Omega\times\mathbb{R}$ be a compact $\tau$-invariant set. Since $(\Omega,\sigma)$ is minimal, $\mathcal K$ projects onto $\Omega$, that is, the continuous map $\pi\colon\mathcal{K}\rightarrow\Omega$, $(\omega,x)\mapsto\omega$ is surjective. In addition,
\begin{equation*}
\mathcal{K}\subseteq \bigcup_{\omega\in\Omega} \big(\{\omega\}\times [\alpha_\mathcal{K}(\omega),\beta_\mathcal{K}(\omega)]\big)\,,
\end{equation*}
where $\alpha_\mathcal{K}(\omega)=\inf\{x\in\mathbb{R}\colon\, (\omega,x)\in\mathcal{K}\}$ and $\beta_\mathcal{K}(\omega)=\sup\{x\in\mathbb{R}\colon\,(\omega,x)\in\mathcal{K}\}$ are, respectively, lower and upper semicontinuous $\tau$-equilibria whose graphs are contained in $\mathcal{K}$. In particular, the residual sets of their continuity points are $\sigma$-invariant. They will be called the \emph{lower and upper delimiter equilibria} of $\mathcal{K}$. The compact $\tau$-invariant set $\mathcal{K}$ is said to be \emph{pinched} if there exists $\omega\in\Omega$ such that the section $(\mathcal{K})_\omega=\{x\colon\,(\omega,x)\in\mathcal{K}\}$ is a singleton. A $\tau$-minimal set $\mathcal{M}\subset\Omega\times\mathbb{R}$ is said to be \emph{hyperbolic attractive} (resp. \emph{repulsive}) if it is uniformly exponentially asymptotically stable at $\infty$ (resp. $-\infty$). Otherwise, it is said to be \emph{nonhyperbolic}.

A function $f\colon\Omega\times\mathbb{R}\rightarrow\mathbb{R}$ is said to be \emph{coercive} ($\mathrm{(Co)}$ for short) if
\begin{equation*}
\lim_{|x|\rightarrow\infty}\frac{f(\omega,x)}{x}=-\infty
\end{equation*}
uniformly on $\Omega$. A stronger definition of coercivity will be needed in part of Section~\ref{sec:6generalframework} and in Section~\ref{sec:7anotherbifurcation}: $f\colon\Omega\times\mathbb{R}\rightarrow\mathbb{R}$ is said to be \emph{$2$-coercive} ($\mathrm{(Co)}_2$ for short) if
\begin{equation*}
\lim_{x\rightarrow\pm\infty}\frac{f(\omega,x)}{x^2}=\mp\infty\,,
\end{equation*}
uniformly on $\Omega$. It is clear that, if $f$ is $\mathrm{(Co)}_2$, then it is $\mathrm{(Co)}$. The arguments leading to Theorem~16 of \cite{pullbackforwards} (see also Section~1.2 of \cite{carvalho1}) ensure that, if $f\in C^{0,1}(\Omega\times\mathbb{R},\mathbb{R})$ is $\mathrm{(Co)}$, then the flow $\tau$ is globally forward defined and admits a global attractor. That is, a compact $\tau$-invariant set $\mathcal{A}$ which satisfies $\lim_{t\rightarrow\infty} \text{dist}(\tau_t(\mathcal{C}),\mathcal{A})=0$ for every bounded set $\mathcal{C}\subset\Omega\times\mathbb{R}$, where $\tau_t(\mathcal{C})=\{(\omega{\cdot}t,u(t,\omega,x_0))\colon\, (\omega,x_0)\in\mathcal{C}\}$ and
\begin{equation*}
\text{dist}(\mathcal{C}_1,\mathcal{C}_2)=\sup_{(\omega_1,x_1)\in\mathcal{C}_1}\left(\inf_{(\omega_2,x_2)\in\mathcal{C}_2}\big(\mathrm{dist}_{\Omega\times\mathbb{R}}((\omega_1,x_1),(\omega_2,x_2))\big)\right)\,.
\end{equation*}
In addition, the attractor takes the form $\mathcal{A}=\bigcup_{\omega\in\Omega}(\{\omega\}\times[\alpha_\mathcal{A}(\omega),\beta_\mathcal{A}(\omega)])$ and is composed by the union of all the globally defined and bounded $\tau$-orbits. And, as proved in Theorem~5.1(iii) of \cite{dno1}, any global (strict) lower solution $\kappa$ satisfies $\kappa\leq \beta_\mathcal{A}$ ($\kappa<\beta_\mathcal{A}$) and a (strict) upper solution $\kappa$ satisfies $\kappa\geq \alpha_\mathcal{A}$ ($\kappa>\alpha_\mathcal{A}$). Besides, if $\omega_1$ (resp.~$\omega_2$) is a continuity point of $\alpha_{\mathcal A}$ (resp.~$\beta_{\mathcal A}$), then the set $\mathcal M^l=\mathrm{cl}_{\Omega\times\mathbb R}\{(\omega_1{\cdot}t,\alpha_{\mathcal A}(\omega_1{\cdot}t))\colon\,t\in\mathbb R\}$ (resp.~$\mathcal M^u=\text{\rm cl}_{\Omega\times\mathbb R}\{(\omega_2{\cdot}t,\beta_{\mathcal A}(\omega_2{\cdot}t))\colon\,t\in\mathbb R\}$) is the lower (resp.~upper) $\tau$-minimal set, and its sections reduce to the points $\alpha_{\mathcal A}(\omega)$ (resp.~$\beta_{\mathcal A}(\omega)$) at all the continuity points $\omega$ of $\alpha_{\mathcal A}$ (resp.~$\beta_\mathcal{A}$): see Theorem 3.3 of \cite{lineardissipativescalar}. Moreover, it is easy to check by contradiction that, if $\mathcal M^l$ (resp.~$\mathcal M^u$) is hyperbolic, then it is attractive and it coincides with the graph of $\alpha_{\mathcal A}$ (resp.~$\beta_{\mathcal A}$), which therefore is a continuous map.

\subsection{Measures and Lyapunov exponents}\label{subsec:2measuresandlyapunovexponents}
We shall say that a Borel measure $m$ on $\Omega$ is \emph{normalized} if $m(\Omega)=1$, that it is $\sigma$-\emph{invariant} if $m(\sigma_t(\mathcal{B}))=m(\mathcal{B})$ for every $t\in\mathbb{R}$ and every Borel subset $\mathcal{B}\subseteq\Omega$, and that it is $\sigma$-\emph{ergodic} if it is normalized, $\sigma$-invariant and $m(\mathcal{B})\in\{0,1\}$ for every $\sigma$-invariant subset $\mathcal{B}\subseteq\Omega$. $\mathfrak{M}_\mathrm{inv}(\Omega,\sigma)$ and $\mathfrak{M}_\mathrm{erg}(\Omega,\sigma)$ are, respectively, the nonempty sets of normalized $\sigma$-invariant and $\sigma$-ergodic Borel measures on $\Omega$. The flow $(\Omega,\sigma)$ is said to be \emph{uniquely ergodic} if $\mathfrak{M}_\text{inv}(\Omega,\sigma)$ reduces to just one element $m$, in which case $m$ is ergodic; and it is said to be \emph{finitely ergodic} if $\mathfrak{M}_\mathrm{erg}(\Omega,\sigma)$ is a finite set.

The \emph{Lyapunov exponent} of $a\in C(\Omega)$ with respect to $m\in\mathfrak{M}_\mathrm{erg}(\Omega,\sigma)$ is
\begin{equation*}
\gamma_a(\Omega,m)=\int_\Omega a(\omega)\, dm\,.
\end{equation*}
The family of scalar linear differential equations $x'=a(\omega{\cdot}t)\,x$ has \emph{exponential dichotomy} over $\Omega$ if there exist $k\geq1$ and $\delta>0$ such that either
\begin{equation*}
\exp\int_0^t a(\omega{\cdot}s)\; ds\leq ke^{-\delta t}\text{ whenever }\omega\in\Omega\text{ and }t\geq0
\end{equation*}
or
\begin{equation*}
\exp\int_0^t a(\omega{\cdot}s)\; ds\leq ke^{\delta t}\text{ whenever }\omega\in\Omega\text{ and }t\leq0\,.
\end{equation*}
The set of $\lambda\in\mathbb{R}$ such that the family $x'=(a(\omega{\cdot}t)-\lambda)\,x$ does not have exponential dichotomy over $\Omega$ is called the \emph{Sacker and Sell spectrum of} $a\in C(\Omega)$, and represented by $\mathrm{sp}(a)$. Recall that $\Omega$ is connected, since $(\Omega,\sigma)$ is minimal. The arguments in \cite{johnson2} and \cite{sackersell3} show the existence of $m^l,m^u\in\mathfrak{M}_\mathrm{erg}(\Omega,\sigma)$ such that $\mathrm{sp}(a)=\left[\gamma_a(\Omega,m^l),\gamma_a(\Omega,m^u)\right]$, and also that $\int_\Omega a(\omega)\,dm\in\mathrm{sp}(a)$ for any $m\in \mathfrak{M}_\mathrm{inv}(\Omega,\sigma)$. We shall say that $a$ has \emph{band spectrum} if $\mathrm{sp}(a)$ is a nondegenerate interval and that $a$ has \emph{point spectrum} if $\mathrm{sp}(a)$ reduces to a point. As seen in Subsection~\ref{subsec:2boundedprimitive}, $\mathrm{sp}(a)=\{0\}$ if $a\in C_0(\Omega)$.

On the other hand, assume that $f\in C^{0,1}(\Omega\times\mathbb{R},\mathbb{R})$, where $f$ is the function on the right hand side of \eqref{eq:generalequation}. The \emph{Lyapunov exponent} of a compact $\tau$-invariant set $\mathcal{K}\subset\Omega\times\mathbb{R}$ with respect to $\nu\in\mathfrak{M}_\mathrm{erg}(\mathcal{K},\tau)$ is
\begin{equation*}
\gamma_{f_x}(\mathcal{K},\nu)=\int_\mathcal{K} f_x(\omega,x)\, d\nu\,.
\end{equation*}
We will frequently omit the subscript $f_x$ if no confusion may arise. We will refer to the Sacker and Sell spectrum of $f_x\colon\mathcal{K}\rightarrow\mathbb{R}$ as the \emph{Sacker and Sell spectrum of $f_x$ on a compact $\tau$-invariant set $\mathcal{K}\subset\Omega\times\mathbb{R}$}. 
Since $(\Omega,\sigma)$ is a minimal flow, a $\tau$-minimal set $\mathcal{M}\subset\Omega\times\mathbb{R}$ is nonhyperbolic if and only if $0$ belongs to the Sacker and Sell spectrum of $f_x$ on $\mathcal{M}$. Moreover, Proposition~2.8 of \cite{lineardissipativescalar} proves that $\mathcal{M}$ is an attractive (resp. repulsive) hyperbolic copy of the base if and only if all its Lyapunov exponents are strictly negative (resp.~positive).

Theorems~1.8.4 of \cite{arnol} and 4.1 of \cite{furstenberg1} provide a fundamental characterization of the set $\mathfrak{M}_\mathrm{erg}(\mathcal{K},\tau)$ given by the $\tau$-ergodic measures concentrated on a compact $\tau$-invariant set $\mathcal K\subset\Omega\times\mathbb{R}$: for any $\nu\in\mathfrak{M}_\mathrm{erg}(\mathcal{K},\tau)$, there exists an $m$-measurable $\tau$-equilibrium $\beta\colon\Omega\rightarrow\mathbb{R}$ with graph contained in $\mathcal{K}$ such that, for every continuous function $g\colon\Omega\times\mathbb{R}\rightarrow\mathbb{R}$,
\begin{equation}\label{eq:new}
\int_\mathcal{K} g(\omega,x)\; d\nu=\int_\Omega g(\omega,\beta(\omega))\; dm\,,
\end{equation}
where $m\in\mathfrak{M}_\mathrm{erg}(\Omega,\sigma)$ is the ergodic measure on which $\nu$ projects, given by $m(A)=\nu((A\times\mathbb{R})\cap \mathcal{K})$. In particular, the Lyapunov exponent on $\mathcal{K}$ for \eqref{eq:generalequation} with respect to any $\tau$-ergodic measure projecting onto $m$ is given by an integral of the form $\int_\Omega f_x(\omega,\beta(\omega))\, dm$. The converse also holds: any $m$-measurable $\tau$-equilibrium $\beta:\Omega\rightarrow\mathbb{R}$ with graph in $\mathcal{K}$ defines $\nu\in\mathfrak{M}_\text{erg}(\mathcal{K},\tau)$ projecting onto $m$ by \eqref{eq:new}. Note that $\beta_1$ and $\beta_2$ define the same measure if and only if they coincide $m$-a.e.
\subsection{Strict d-concavity}\label{subsec2:dconcavity} We shall say that $f\in C^{0,1}(\Omega\times\mathbb{R},\mathbb{R})$ is \emph{d-concave} ($\mathrm{(DC)}$ for short) if its derivative $f_x$ is concave on $\mathbb{R}$ for all $\omega\in\Omega$. With the purpose of measuring the degree of strictness of the concavity of $f_x$, the standardized $\epsilon$-modules of d-concavity of $f$ on a compact interval $J$ were introduced in \cite{dno1}, and several subsets of strictly d-concave functions of $C^{0,1}(\Omega\times\mathbb{R},\mathbb{R})$ were defined in terms of these modules. In this paper, we will only be interested in the set $\mathrm{(SDC)_*}$ of \emph{strictly d-concave functions with respect to every measure} (see Definition~3.8 of \cite{dno1}). Proposition~3.9 of \cite{dno1} gives a characterization of this set of functions which will be sufficient for the purposes of this paper: $f\in C^{0,2}(\Omega\times\mathbb{R},\mathbb{R})$ is $\mathrm{(SDC)}_*$ if and only if $m(\{\omega\in\Omega\colon\, f_{xx}(\omega,\cdot)$ is strictly decreasing on $J\})>0$ for every compact interval $J$ and every $m\in\mathfrak{M}_\mathrm{erg}(\Omega,\sigma)$. In particular, it can be easily checked that any polynomial of the form $p(\omega,x)=-a_3(\omega)x^3+a_2(\omega)x^2+a_1(\omega)x+a_0(\omega)$, where the coefficients are continuous and $\Omega$ and $a_3$ is nonnegative and nonzero, is $\mathrm{(SDC)}_*$, since $p_{xx}(\omega,\cdot)$ is strictly decreasing on $\mathbb{R}$ for every $\omega$ on an open subset of $\Omega$ (recall that the minimality of $(\Omega,\sigma)$ ensures that every open set has positive $m$-measure for all $m\in\mathfrak{M}_\mathrm{erg}(\Omega,\sigma)$).

Assume that the function $f$ of \eqref{eq:generalequation} is $\mathrm{(SDC)}_*$. Following the methods of \cite{pliss1} and \cite{tineo1}, Theorems~4.1 and 4.2 of \cite{dno1} state relevant dynamical properties of the local skewproduct flow $\tau$ in terms of the previous properties. Let $\mathcal{K}\subset\Omega\times\mathbb{R}$ be a compact $\tau$-invariant set. Then, there exist at most three distinct $\tau$-invariant measures of $\mathfrak{M}_\mathrm{erg}(\mathcal{K},\tau)$ which project onto $m$. Moreover, if there exist three such measures $\nu_1$, $\nu_2$ and $\nu_3$ projecting onto $m$, and they are respectively given by the $m$-measurable equilibria $\beta_1$, $\beta_2$ and $\beta_3$ (see \eqref{eq:new}) with $\beta_1(\omega)<\beta_2(\omega)<\beta_3(\omega)$ for $m$-a.e. $\omega\in\Omega$, then $\gamma_{f_x}(\mathcal{K},\nu_1)<0$, $\gamma_{f_x}(\mathcal{K},\nu_2)>0$ and $\gamma_{f_x}(\mathcal{K},\nu_3)<0$ (see the proof of Theorem~4.1 of \cite{dno1}). In addition, $\mathcal{K}$ contains at most three disjoint compact $\tau$-invariant sets, and if it contains exactly three, then they are hyperbolic copies of the base: attractive the upper and lower ones, and repulsive the middle one. These properties will be often combined with those established in Proposition~5.3 of \cite{dno1}: if $f$ is coercive and either if there exists a repulsive hyperbolic $\tau$-minimal set or if there exist two hyperbolic $\tau$-minimal sets, then there exist three $\tau$-minimal sets.
\section{Generalized pitchfork bifurcation patterns}\label{sec:3bifurcationtheorem}
Let $(\Omega,\sigma)$ be a minimal flow, and let $f\in C^{0,2}(\Omega\times\mathbb{R},\mathbb{R})$ be $\mathrm{(Co)}$ and $\mathrm{(SDC)}_*$, and
satisfy $f(\omega,0)=0$ for every $\omega\in\Omega$. The description of the global bifurcation diagram for the family
\begin{equation}\label{eq:65}
x'=f(\omega{\cdot}t,x)+\lambda x\,,\quad\omega\in\Omega
\end{equation}
with respect the real parameter $\lambda$ was initiated in Section~6 of \cite{dno1} and is completed in Theorem \ref{th:clasification}. We denote by $\mathcal{I}_{\omega,x_0}^\lambda\rightarrow\mathbb{R}$, $t\mapsto u_\lambda(t,\omega,x_0)$ the maximal solution of $\eqref{eq:65}_\lambda$ with initial value $u_\lambda(0,\omega,x_0)=x_0$, and by $\tau_\lambda$ the corresponding local skewproduct flow induced by $\eqref{eq:65}_\lambda$; i.e., $\tau_\lambda(t,\omega,x_0)=(\omega{\cdot}t,u_\lambda(t,\omega,x_0))$. The coercivity property ensures the existence of the global attractor $\mathcal{A}_\lambda$ for all $\lambda\in\mathbb{R}$. Theorems~5.1 and 6.2 of \cite{dno1} describe the structure of $\mathcal{A}_\lambda$ and its variation with respect to $\lambda$: it can be written as
\begin{equation*}
\mathcal{A}_\lambda=\bigcup_{\omega\in\Omega} \big(\{\omega\}\times[\alpha_\lambda(\omega),\beta_\lambda(\omega)]\big)\,,
\end{equation*}
where: $\alpha_\lambda\colon\Omega\to\mathbb R$ and $\beta_\lambda\colon\Omega\to\mathbb R$ are lower and upper semicontinuous $\tau_\lambda$-equilibria;
$\alpha_\lambda(\omega)\leq0\leq\beta_\lambda(\omega)$ for all $\omega\in\Omega$; the maps $\lambda\mapsto\beta_\lambda(\omega)$ and $\lambda\mapsto\alpha_\lambda(\omega)$ are respectively nondecreasing and nonincreasing on $\mathbb{R}$ for all $\omega\in\Omega$ and both are right-continuous; $\lim_{\lambda\rightarrow\infty}\alpha_\lambda(\omega)=-\infty$ and $\lim_{\lambda\rightarrow\infty}\beta_\lambda(\omega)=\infty$ uniformly on $\Omega$; and there exists $\lambda_*\in\mathbb{R}$ such that $\mathcal{A}_\lambda=\Omega\times\{0\}$ and it is hyperbolic attractive for every $\lambda<\lambda_*$. We represent by $\mathcal M^l_\lambda$ and $\mathcal M^u_\lambda$ the lower and upper $\tau_\lambda$-minimal sets, and recall that any of them is hyperbolic if and only it is hyperbolic attractive, in which case it coincides with the graph of the corresponding delimiter equilibrium of $\mathcal A_\lambda$ (see Subsection~\ref{subsec:2compactcoerciveattractor}).

Let $\{\kappa_\lambda^1:\lambda\in(\lambda_1,\lambda_2)\}$ and $\{\kappa_\lambda^2:\lambda\in(\lambda_1,\lambda_2)\}$ be two families of $\tau_\lambda$-equilbria. We will say that $\kappa^1_\lambda$ and $\kappa^2_\lambda$ {\em collide} (or that {\em $\kappa^1_\lambda$ collides with} $\kappa^2_\lambda$) {\em as $\lambda\downarrow\lambda_1$ or $\lambda\uparrow\lambda_2$ on a residual set $\mathcal{R}\subseteq\Omega$} if the limit of the difference $\kappa_\lambda^1(\omega)-\kappa_\lambda^2(\omega)$ vanishes for all $\omega\in\mathcal{R}$. The same terminology will be used for the different parametric families of flows appearing throughout the paper.
\begin{theorem}\label{th:clasification} Let $f\in C^{0,2}(\Omega\times\mathbb{R},\mathbb{R})$ be $\mathrm{(Co)}$ and $\mathrm{(SDC)}_*$, and let $[-\lambda_+,-\lambda_-]$ be the Sacker and Sell spectrum of $f_x$ on $\mathcal{M}_0=\Omega\times\{0\}$, with $\lambda_-\leq\lambda_+$. Then, $\mathcal{M}_0$ is hyperbolic attractive (resp. repulsive) if $\lambda<\lambda_-$ (resp. $\lambda>\lambda_+$) and nonhyperbolic if $\lambda\in[\lambda_-,\lambda_+]$; $\tau_\lambda$ admits three different hyperbolic minimal sets $\mathcal{M}_\lambda^l<\mathcal{M}_0<\mathcal{M}_\lambda^u$ for $\lambda>\lambda_+$, where $\mathcal{M}_\lambda^l$ and $\mathcal{M}_\lambda^u$ are hyperbolic attractive and given by the graphs of $\alpha_\lambda$ and $\beta_\lambda$ respectively; and either $\alpha_\lambda$ or $\beta_\lambda$ (or both of them) collide with $0$ on a residual $\sigma$-invariant set as $\lambda\downarrow\lambda_+$.

Assume that this is the case for $\beta_\lambda$. Then, one of the following situations holds:
\begin{enumerate}[label=\rm{(\roman*)}]
\item {\rm(Local saddle-node and transcritical bifurcations).} $\lambda_-\leq\lambda_+$ and there exists $\lambda_0<\lambda_-$ such that: $\mathcal{A}_\lambda=\mathcal{M}_0$ for $\lambda<\lambda_0$; $\tau_{\lambda_0}$ admits exactly two different minimal sets $\mathcal{M}_{\lambda_0}^l<\mathcal{M}_0$, with $\mathcal{M}_{\lambda_0}^l$ nonhyperbolic; $\tau_\lambda$ admits three hyperbolic minimal sets $\mathcal{M}_\lambda^l<\mathcal{N}_\lambda<\mathcal{M}_0$ for $\lambda\in(\lambda_0,\lambda_-)$, where $\mathcal{M}_\lambda^l$ is hyperbolic attractive and given by the graph of $\alpha_\lambda$, and $\mathcal{N}_\lambda$ is hyperbolic repulsive and given by the graph of a continuous map $\kappa_\lambda\colon\Omega\rightarrow\mathbb{R}$ which increases strictly as $\lambda$ increases in $(\lambda_0,\lambda_-)$, and which collides with $\alpha_\lambda$ (resp. with $0$) on a residual $\sigma$-invariant set as $\lambda\downarrow\lambda_0$ (resp. $\lambda\uparrow\lambda_-$); and $\tau_\lambda$ admits exactly two minimal sets $\mathcal{M}_\lambda^l<\mathcal{M}_0$ for $\lambda\in[\lambda_-,\lambda_+]$, where $\mathcal{M}_\lambda^l$ is hyperbolic attractive and given by the graph of $\alpha_\lambda$. In particular, $\lambda_0$, $\lambda_-$ and $\lambda_+$ are the unique bifurcation points: a local saddle-node bifurcation of minimal sets occurs around $\mathcal{M}_{\lambda_0}$ at $\lambda_0$, as well as a discontinuous bifurcation of attractors; and a classical (resp. generalized) transcritical bifurcation of minimal sets arises around $\mathcal{M}_0$ at $\lambda_-$ (resp. on $[\lambda_-,\lambda_+]$) if $\lambda_-=\lambda_+$ (resp. $\lambda_-<\lambda_+$).
\item {\rm(Global classical pitchfork bifurcation).} $\lambda_-\leq\lambda_+$, and $\mathcal{M}_0$ is the unique $\tau_\lambda$-minimal set if $\lambda\leq\lambda_+$, and $\mathcal{A}_\lambda=\mathcal{M}_0$ if $\lambda<\lambda_-$; and both $\alpha_\lambda$ and $\beta_\lambda$ collide with $0$ on a residual $\sigma$-invariant set as $\lambda\downarrow\lambda_+$. A classical pitchfork bifurcation of minimal sets arises around $\mathcal{M}_0$ at $\lambda_+$.
\item {\rm(Global generalized pitchfork bifurcation).} $\lambda_-<\lambda_+$ and there exists $\lambda_0\in[\lambda_-,\lambda_+)$ such that: $\mathcal{M}_0$ is the unique $\tau_\lambda$-minimal set if $\lambda<\lambda_0$ and $\mathcal{A}_\lambda=\mathcal{M}_0$ if $\lambda<\lambda_-$; and $\tau_\lambda$-admits two minimal sets $\mathcal{M}_\lambda^l<\mathcal{M}_0$ for $\lambda\in(\lambda_0,\lambda_+]$, where $\mathcal{M}_\lambda^l$ is hyperbolic attractive and given by the graph of $\alpha_\lambda$. A generalized pitchfork bifurcation of minimal sets around $\mathcal{M}_0$ on $[\lambda_-,\lambda_+]$ arises, with $\lambda_0$ and $\lambda_+$ as bifurcation points.
\end{enumerate}
The possibilities for the global bifurcation diagram are symmetric with respect to the horizontal axis to those described if $\alpha_\lambda$ collides with $0$ as $\lambda\downarrow\lambda_+$.
\end{theorem}
The following technical lemma is needed in our proof of the theorem.
\begin{lemma} \label{lemma:tapaslambdax} Let $f\in C^{0,2}(\Omega\times\mathbb{R},\mathbb{R})$ be $\mathrm{(Co)}$ and $(\mathrm{SDC})_*$ and satisfy $f(\omega,0)=0$ for all $\omega\in\Omega$. Assume that there exist exactly two minimal sets $\mathcal M^l<\mathcal M_0$ (resp. $\mathcal M^u>\mathcal M_0$) for the local skewproduct flow $\tau$ defined by the solutions $u(t,\omega,x_0)$ of the family $x'=f(\omega{\cdot}t,x)$, with $\mathcal{M}^l$ (resp. $\mathcal{M}^u$) hyperbolic attractive. Assume also that $\int_\Omega f_x(\omega,0)\,dm\ne 0$ for an $m\in\mathfrak{M}_\mathrm{erg}(\Omega,\sigma)$. Then, there exist a lower semicontinuous $\tau$-equilibrium $\kappa_1\colon\Omega\rightarrow(-\infty,0]$ and an upper semicontinuous $\tau$-equilibrium $\kappa_2\colon\Omega\rightarrow[0,\infty)$ which take the value $0$ on the residual $\sigma$-invariant sets of their continuity points and such that $\int_\Omega f_x(\omega,\kappa_1(\omega))\,dm>0$ and $\int_\Omega f_x(\omega,\kappa_2(\omega))\,dm<0$ (resp.~$\int_\Omega f_x(\omega,\kappa_1(\omega))\,dm<0$ and $\int_\Omega f_x(\omega,\kappa_2(\omega))\,dm>0$). Moreover,
\begin{equation*}
\begin{split}
\kappa_1(\omega)&=\sup\{x_0\in\mathbb{R}\colon\,\lim_{t\rightarrow\infty}(u(t,\omega,x_0)-\alpha(\omega{\cdot}t))=0\}>\alpha(\omega)\\
(\text{resp.}\;\kappa_2(\omega)&=\inf\{x_0\in\mathbb{R}\colon\,\lim_{t\rightarrow\infty}(u(t,\omega,x_0)-\beta(\omega{\cdot}t))=0\}<\beta(\omega))\,,\\
\end{split}
\end{equation*}
where the $\tau$-equilibrium $\alpha$ (resp. $\beta$) is the lower (resp. upper) delimiter of the global attractor, with graph $\mathcal{M}^l$ (resp. $\mathcal{M}^u$).
\end{lemma}
\begin{proof} Proposition 5.3(ii) of \cite{dno1} ensures the nonhyperbolicity of $\mathcal{M}_0$. We reason in the case $\mathcal{M}^l<\mathcal{M}_0$. Theorem~5.12 of \cite{dno1} ensures that the global $\xi$-bifurcation diagram for $x'=f(\omega{\cdot}t,x)+\xi$ is described by Theorem~5.10 of \cite{dno1}, with $0$ as left local saddle-node bifurcation point, at which the upper (nonhyperbolic) minimal set is $\mathcal{M}_0$. As explained in the proof of Theorem~5.10 of \cite{dno1}, $\mathcal{M}_0$ is contained on a compact $\tau$-invariant pinched set $\bigcup_{\omega\in\Omega} \big(\{\omega\}\times [\kappa_1(\omega),\kappa_2(\omega)]\big)$, where: $\kappa_1\colon\Omega\rightarrow(-\infty,0]$ (called $\kappa$ in \cite{dno1}) is a lower semicontinuous $\tau$-equilibrium; $\kappa_2\colon\Omega\rightarrow(-\infty,0]$ (called $\beta$ in \cite{dno1}) is the upper delimiter of the global attractor; and both of them coincide (and hence take the value 0) on the residual set of their common continuity points (and hence at all their continuity points: see Proposition 2.4 of \cite{dno1}). Since $\int_\Omega f_x(\omega,0)\, dm\ne 0$, Proposition~5.11(i) of \cite{dno1} ensures that $\int_\Omega f_x(\omega,\kappa_1(\omega))\,dm>0$ and $\int_\Omega f_x(\omega,\kappa_2(\omega))\,dm\,<0$.

Let us check the last assertion, also in the case $\mathcal{M}^l<\mathcal{M}_0$. Recall that the hyperbolic attractiveness of $\mathcal{M}^l$ ensures that it coincides with the graph of $\alpha$, which is continuous, and that every global strict upper solution is strictly greater than $\alpha$ (see Subsection~\ref{subsec:2compactcoerciveattractor}). Since $\kappa_1$ is a $\tau$-equilibrium, $\kappa_1(\omega)\geq\sup\{x_0\in\mathbb{R}\colon\,\lim_{t\rightarrow\infty}(u(t,\omega,x_0)-\alpha(\omega{\cdot}t))=0\}$ for all $\omega\in\Omega$. We fix $\omega_0$ to check that $\lim_{t\rightarrow\infty}(u(t,\omega_0,x_0)-\alpha(\omega_0{\cdot}t))=0$ for any $x_0<\kappa_1(\omega_0)$. The global $\xi$-bifurcation diagram for $x'=f(\omega{\cdot}t,x)+\xi$ described by Theorem~5.10 of \cite{dno1} ensures the existence of $\xi_+>0$ and strictly negative continuous equilibria $\kappa_\xi\colon\Omega\rightarrow\mathbb{R}$ of $x'=f(\omega{\cdot}t,x)+\xi$ for $\xi\in(0,\xi_+)$ which decrease strictly as $\xi$ increases and such that $\kappa_1(\omega)=\lim_{\xi\downarrow0}\kappa_\xi(\omega)$ for all $\omega\in\Omega$. Hence, for any $x_0<\kappa_1(\omega_0)$, there exists $\xi_0>0$ such that $x_0<\kappa_{\xi_0}(\omega_0)$. Since $\kappa_{\xi_0}$ is a strict upper solution of $x'=f(\omega{\cdot}t,x)$ and hence a strong $\tau$-superequilibrium, the $\boldsymbol\upomega$-limit for $\tau$ of $(\omega_0,x_0)$ is strictly below the graph of $\kappa_{\xi_0}$ (see Subsection~\ref{subsec:2equilibria}), so this $\boldsymbol\upomega$-limit contains $\mathcal{M}^l$, which is the unique $\tau$-minimal set strictly below it. Therefore, there exists a sequence $\{t_n\}\uparrow\infty$ such that $\lim_{n\rightarrow\infty}(u(t_n,\omega_0,x_0)-\alpha(\omega_0{\cdot}t_n))=0$, and the hyperbolic attractiveness of $\mathcal{M}^l$ ensures that $\lim_{t\rightarrow\infty}(u(t,\omega_0,x_0)-\alpha(\omega_0{\cdot}t))=0$, as asserted.

The proof is analogous in the other case.
\end{proof}
\begin{proof}[Proof of Theorem~$\mathrm{\ref{th:clasification}}$] The Sacker and Sell spectrum of $f_x+\lambda$ on $\mathcal M_0$ is $[-\lambda_++\lambda, -\lambda_-+\lambda]$, which ensures the stated hyperbolicity properties of $\mathcal M_0$ (see Subsection \ref{subsec:2measuresandlyapunovexponents}). As in Theorem~6.3 of \cite{dno1}, we define
\begin{equation*}
\begin{split}
\mu_-=\inf\{\lambda\colon\forall\;\xi>\lambda \text{ the graph of }\alpha_\xi\text{ is the hyperbolic minimal set }\mathcal{M}_\xi^l<\mathcal{M}_0\}\,,\\
\mu_+=\inf\{\lambda\colon\forall\;\xi>\lambda \text{ the graph of }\beta_\xi\text{ is the hyperbolic minimal set }\mathcal{M}_\xi^u>\mathcal{M}_0\}\,,\\
\end{split}
\end{equation*}
which, as proved there, belong to $(-\infty,\lambda_+]$. This property guarantees the stated structure of the $\tau_\lambda$-minimal sets for $\lambda>\lambda_+$, since there exist at most three $\tau_\lambda$-minimal sets (see Subsection~\ref{subsec2:dconcavity}). Theorem~6.3(ii) of \cite{dno1} also ensures that at least one of these two parameters $\mu_-$, $\mu_+$ coincides with $\lambda_+$, which proves the stated collision properties for $\alpha_\lambda$ or for $\beta_\lambda$ as $\lambda\downarrow\lambda_+$. As in the statement, we assume that this is the case for $\beta_\lambda$, i.e., that $\mu_+=\lambda_+$. Then, since $\lambda\mapsto\beta_\lambda(\omega)$ is nondecreasing for all $\omega\in\Omega$ and the intersection of two residual sets is also residual, $\mathcal{M}_0$ is the upper minimal set for all $\lambda\leq\lambda_+$. If also $\mu_-=\lambda_+$, then Theorem~6.3(iii) of \cite{dno1} ensures that the bifurcation diagram is that of (ii). If $\mu_-<\lambda_-$, then Theorem~6.4 of \cite{dno1} shows that the diagram is that of (i), with $\lambda_0=\mu_-$. The remaining case is, hence, $\mu_-\in[\lambda_-,\lambda_+)$. We will check that, in this case, the situation is that of (iii), which will complete the proof.

Let us call $\lambda_0=\mu_-$. Notice that, if $\lambda\in(\lambda_0,\lambda_+]$, then there exist only two $\tau_\lambda$-minimal sets, as otherwise the nonhyperbolicity of $\mathcal{M}_0$ would be contradicted (see Subsection~\ref{subsec2:dconcavity}); and, as explained before, $\mathcal M_\lambda^l$ is hyperbolic attractive (and given by the graph of $\alpha_\lambda$). Consequently, it only remains to prove that $\alpha_\lambda(\omega)=0$ on the residual $\sigma$-invariant set of its continuity points for $\lambda<\lambda_0$. This will ensure that $\mathcal{M}_0$ is the unique $\tau_\lambda$-minimal set for $\lambda<\lambda_0$, and the hyperbolic attractiveness of $\mathcal{M}_0$ for $\lambda<\lambda_-$ will ensure that $\mathcal{A}_\lambda=\mathcal{M}_0$ for $\lambda<\lambda_-$ (see Theorem~3.4 of \cite{lineardissipativescalar}). Recall also that
$\alpha_\lambda$ vanishes at all its continuity points if it vanishes at one of them (see e.g.~Proposition~2.5 of \cite{dno1}).

First, let us assume that $\mathcal{M}_{\lambda_0}^l=\mathcal{M}_0$, which means that $\alpha_{\lambda_0}(\omega)=0$ on the residual set of its continuity points (see Subsection \ref{subsec:2compactcoerciveattractor}). Therefore, the same happens with $\alpha_\lambda$ if $\lambda<\lambda_0$, since $\alpha_{\lambda_0}\leq\alpha_\lambda\leq0$ and the intersection of two residual sets is also residual. This proves the result in this case.

So, we assume $\mathcal{M}_{\lambda_0}^l<\mathcal{M}_0$, which in particular ensures that $\alpha_{\lambda_0}(\omega)<0$ for all $\omega\in\Omega$. Let us take $\lambda>\lambda_0$. Lemma~\ref{lemma:tapaslambdax} shows that the upper delimiter of the basin of attraction of the graph of $\alpha_\lambda$ is a lower semicontinuous $\tau_\lambda$-equilibrium $\kappa^1_\lambda\colon\Omega\rightarrow(-\infty,0]$ which vanishes at its continuity points, with $\kappa^1_\lambda>\alpha_\lambda$, and which satisfies $\int_\Omega (f_x(\omega,\kappa^1_\lambda(\omega))+\lambda)\, dm>0$ whenever $\int_\Omega (f_x(\omega,0)+\lambda)\, dm<0$. In particular, since $m$ is ergodic, $\kappa^1_\lambda<0$ $m$-a.e. whenever $\int_\Omega (f_x(\omega,0)+\lambda)\, dm<0$. Let us check that the map $\lambda\mapsto\kappa^1_\lambda(\omega)$ is nondecreasing on $(\lambda_0,\infty)$ for all $\omega\in\Omega$. It is easy to deduce from the hyperbolicity of the graph $\mathcal M_{\lambda}^l$ of $\alpha_\lambda$ that $x<\kappa_{\lambda}^1(\w)$ if and only if $\mathcal M_{\lambda}^l$ is the $\boldsymbol\upomega$-limit for $\tau_\lambda$ of $(\omega,x)$. We take $\lambda_0<\lambda_1<\lambda_2$, $\omega\in\Omega$ and $x<\kappa^1_{\lambda_1}(\omega)$. Then, $u_{\lambda_2}(t,\omega,x)<u_{\lambda_1}(t,\w,x)<\kappa^1_{\lambda_1}(\omega{\cdot}t)$ for all $t\ge 0$, which precludes the possibility that $\mathcal M_0$ is contained in the $\boldsymbol\upomega$-limit set for $\tau_{\lambda_2}$ of $(\omega,x)$ and hence guarantees that $x<\kappa^1_{\lambda_2}(\w)$. This proves the assertion. As a consequence of this nondecreasing character, the map $\kappa^1_{\lambda_0}(\omega)=\lim_{\lambda\downarrow\lambda_0}\kappa^1_\lambda(\omega)$ is well defined, and it satisfies $\alpha_{\lambda_0}(\omega)\le\kappa^1_{\lambda_0}(\omega)\le \kappa^1_{\lambda}(\omega)$ for all $\omega\in\Omega$ if $\lambda>\lambda_0$. It is clear that $\kappa^1_{\lambda_0}$ it is an $m$-measurable $\tau_{\lambda_0}$-equilibrium for every $m\in\mathfrak{M}_\mathrm{erg}(\Omega,\sigma)$.

Let $\beta_{\mathcal{M}}$ be the upper delimiter $\tau_{\lambda_0}$-equilibrium of $\mathcal{M}_{\lambda_0}^l$. Our next purpose is to check that there exist points $\omega_0\in\Omega$ with $\kappa_{\lambda_0}^1(\omega_0)\le\beta_{\mathcal M}(\omega_0)$.
Since $\mathcal{M}_{\lambda_0}^l$ is nonhyperbolic, there exists $m\in\mathfrak{M}_\mathrm{erg}(\Omega,\sigma)$ and an $m$-measurable $\tau_{\lambda_0}$-equilibrium $\tilde\kappa$ with graph in $\mathcal{M}_{\lambda_0}^l$ such that $\int_\Omega(f_x(\omega,\tilde\kappa(\omega))+\lambda_0)\, dm\geq 0$ (see Subsection~\ref{subsec:2measuresandlyapunovexponents}). Proposition~4.4 of \cite{dno1} ensures that $\int_\Omega(f_x(\omega,0)+\lambda_0)\, dm<0$, so $\int_\Omega(f_x(\omega,0)+\lambda)\, dm<0$ for $\lambda\geq\lambda_0$ close enough. Therefore, as seen before, $\kappa^1_\lambda(\omega)<0$ $m$-a.e. for these values of the parameter, and hence $\kappa^1_{\lambda_0}(\omega)<0$ $m$-a.e. Assume for contradiction that $\kappa_{\lambda_1}^1(\omega)>\beta_{\mathcal M}(\omega)$ for all $\omega\in\Omega$, and hence that $\kappa_{\lambda_1}^1(\omega)>\tilde\kappa(\omega)$ for all $\omega\in\Omega$. Then, $0$, $\kappa^1_{\lambda_0}$ and $\tilde\kappa$ define three different $\tau_{\lambda_0}$-ergodic measures by \eqref{eq:new}. This ensures that $\int_\Omega(f_x(\omega,\tilde\kappa(\omega))+\lambda_0)\,dm<0$ (see Subsection~\ref{subsec2:dconcavity}), which is not the case. This proves the assertion.

Finally, let us fix $\lambda<\lambda_0$ and check that $\alpha_\lambda(\omega)=0$ at one of its continuity points, which completes the proof. Propositions~6.1 and 2.1 of \cite{dno1} show the existence of $s^*>0$ such that $\beta_{\mathcal{M}}(\omega)<u_\lambda(s^*,\omega{\cdot}(-s^*),\alpha_{\lambda_0}(\omega{\cdot}(-s^*)))\leq\alpha_\lambda(\omega)$ for all $\omega\in\Omega$ (the last inequality follows from the flow monotonicity and the nonincreasing character of $\lambda\mapsto\alpha_\lambda(\omega)$). We take $\omega_0\in\Omega$ with $\kappa_{\lambda_0}^1(\omega_0)\le\beta_{\mathcal M}(\omega_0)$. Then, $\lim_{\lambda\downarrow\lambda_0}\kappa^1_\lambda(\omega_0)=\kappa^1_{\lambda_0}(\omega_0)\leq\beta_\mathcal{M}(\omega_0)<\alpha_\lambda(\omega_0)$. Therefore, there exists $\xi>\lambda_0$ close enough to get $\kappa^1_\xi(\omega_0)<\alpha_\lambda(\omega_0)$. Consequently, since $\xi>\lambda$ and hence the nonpositive semicontinuous map $\kappa^1_\xi$ is a global lower solution for \eqref{eq:65}$_\lambda$,
$\kappa^1_\xi(\omega_0{\cdot}t)\le u_\lambda(t,\omega_0,\kappa^1_\xi(\omega_0))<u_\lambda(t,\omega_0,\alpha_\lambda(\omega_0))=\alpha_\lambda(\omega_0{\cdot}t)$ for all $t>0$.
Let $\omega_1$ be a common continuity point of $\kappa^1_\xi$ and $\alpha_\lambda$. As $(\Omega,\sigma)$ is minimal, there exists a sequence $\{t_n\}\uparrow\infty$ such that $\omega_0{\cdot}t_n\rightarrow\omega_1$ as $n\rightarrow\infty$. Since $\kappa^1_\xi$ vanishes at $\omega_1$, $0=\kappa^1_\xi(\omega_1)=\lim_{n\rightarrow\infty} \kappa^1_\xi(\omega_0{\cdot}t_n)\leq
\lim_{n\rightarrow\infty}\alpha_\lambda(\omega_0{\cdot}t_n)=\alpha_\lambda(\omega_1)\leq 0$. That is, $\alpha_\lambda(\omega_1)=0$, and the proof is complete.
\end{proof}
There are simple autonomous examples giving rise to situations (ii) (as $x'=-x^3+\lambda x$, with $\lambda_\pm=0$ as bifurcation point, of classical pitchfork type) and (i) (as $x'=-x^3\pm 2x^2+\lambda x$, with $\lambda_0=-1$ as local saddle-node bifurcation point and $\lambda_\pm=0$ as local classical transcritical bifurcation point; the two possibilities of (i) correspond to the two signs of the second-order term). Clearly, case (iii) cannot occur in an autonomous (and hence uniquely ergodic) case. We will go deeper in this matter in Sections~\ref{sec:5criteriacubic}, \ref{sec:6generalframework} and \ref{sec:7anotherbifurcation}, where we will show that all the possibilities realize for suitable families \eqref{eq:65}.
\section{Criteria for cubic polynomial equations}\label{sec:5criteriacubic}
Let us consider families of cubic polynomial ordinary differential equations
\begin{equation}\label{eq:4fcubic}
x'=-a_3(\omega{\cdot}t)x^3+a_2(\omega{\cdot}t)x^2+(a_1(\omega{\cdot}t)+\lambda)x\,,\quad\omega\in\Omega\,,
\end{equation}
where $a_i\in C(\Omega)$ for $i\in\{1,2,3\}$, $a_3$ is strictly positive and $\lambda\in\mathbb{R}$. It is clear that the function $f(\omega,x)=-a_3(\omega)x^3+a_2(\omega)x^2+a_1(\omega)x$ is (Co). In addition, $f$ is $\mathrm{(SDC)}_*$ (see Subsection~\ref{subsec2:dconcavity}). Then, Theorem~\ref{th:clasification} describes the three possible $\lambda$-bifurcation diagrams for \eqref{eq:4fcubic}. Our first goal in this section, achieved in Subsections~\ref{subsec:4CP} and \ref{subsec:4signpreserving}, is to describe conditions on the coefficients $a_i$ determining the specific diagram. The last subsection is devoted to explain how to get actual patterns satisfying the previously established conditions.

As in Section~\ref{sec:3bifurcationtheorem}, $\tau_\lambda(t,\omega,x_0)=(\omega{\cdot}t,u_\lambda(t,\omega,x_0))$ repesents the local skewproduct flow induced by $\eqref{eq:4fcubic}_\lambda$ on $\Omega\times\mathbb{R}$. Notice that the Sacker and Sell spectrum of $f_x$ on $\mathcal{M}_0=\Omega\times\{0\}$ (which is $\tau_\lambda$-minimal for all $\lambda\in\mathbb{R}$) coincides with $\mathrm{sp}(a_1)$.

\subsection{The case of \texorpdfstring{$a_1$}{a1} with continuous primitive}\label{subsec:4CP}
Throughout this subsection, we assume that $a_1\in CP(\Omega)$. Since $CP(\Omega)\subseteq C_0(\Omega)$ (see Subsection~\ref{subsec:2boundedprimitive}), the Sacker and Sell spectrum of $a_1$ is $\mathrm{sp}(a_1)=\{0\}$. Hence, the bifurcation diagram of \eqref{eq:4fcubic} fits in (i) or (ii) of Theorem~\ref{th:clasification}, and our objective is to give criteria ensuring each one of these two possibilities. The relevant fact in terms of which the criteria will be constructed is that the number of $\tau_0$-minimal sets distinguishes the type of bifurcation: there is either one $\tau_0$-minimal set in (ii) or two $\tau_0$-minimal sets in (i).

Proposition~\ref{prop:4Bpextended} provides a simple classification of the casuistic for \eqref{eq:4fcubic} in this case. It is based on the previous bifurcation analysis of
\begin{equation}\label{eq:4cuadratic}
x'=-a_3(\omega{\cdot}t)x^3+(a_2(\omega{\cdot}t)+\xi)x^2\,,\quad\omega\in\Omega\,,
\end{equation}
made in Proposition~\ref{prop:newdiagram}. These two results extend Proposition~6.6 and Corollary~6.7 of \cite{dno1} to the case of strictly positive $a_3$ (instead of $a_3\equiv1$), since the case of $a_1\equiv 0$ is trivially covered by Proposition~\ref{prop:4Bpextended} with $b\equiv0$. We will call $\check\tau_\xi$ the local skewproduct flow induced by $\eqref{eq:4cuadratic}_\xi$ on $\Omega\times\mathbb{R}$.
\begin{proposition}\label{prop:newdiagram} The $\check\tau_\xi$-minimal set $\mathcal{M}_0=\Omega\times\{0\}$ is nonhyperbolic for all $\xi\in\mathbb{R}$. In addition, if $\mathrm{sp}(a_2)=[-\xi_+,-\xi_-]$, with $\xi_-\leq\xi_+$, then
\begin{enumerate}[label=\rm{(\roman*)}]
\item $\check\tau_\xi$ admits exactly two minimal sets $\mathcal{M}_\xi^l<\mathcal{M}_0$ respectively given by the graphs of $\alpha_\xi$ and $0$ for $\xi<\xi_-$, where $\mathcal{M}_\xi^l$ is hyperbolic attractive; and $\alpha_\xi$ collides with $0$ on a residual $\sigma$-invariant set as $\xi\uparrow\xi_-$;
\item $\mathcal{M}_0$ is the unique $\check\tau_\xi$-minimal set for $\xi\in[\xi_-,\xi_+]$;
\item $\check\tau_\xi$ admits exactly two minimal sets $\mathcal{M}_0<\mathcal{M}_\xi^u$ respectively given by the graphs of $0$ and $\beta_\xi$ for $\xi>\xi_+$, where $\mathcal{M}_\xi^u$ is hyperbolic attractive; and $\beta_\xi$ collides with $0$ on a residual $\sigma$-invariant set as $\xi\downarrow\xi_+$.
\end{enumerate}
\end{proposition}
\begin{proof} The proof follows step by step that of Proposition~6.6 of \cite{dno1}. The only remarkable differences arise in checking the existence of a bounded solution of $\eqref{eq:4cuadratic}_\xi$ for $\xi>\xi_+$ which is bounded away from $0$ for $t>0$. Let us explain these differences. Take $\alpha\in(0,\xi-\xi_+)$ with $\alpha<\sqrt{2r}$ for $r=\|a_3\|=\max_{\omega\in\Omega} a_3(\omega)>0$. Let us check that the solution $w(t)$ of $w'=-(a_2(\omega{\cdot}t)+\xi)+a_3(\omega{\cdot}t)/w$ with $w(0)=\alpha<2r/\alpha$ is bounded on $[0,\infty)$. We define $t_1=\sup\{t>0\colon\, w(s)\leq 2r/\alpha+lt_\alpha$ for all $s\in[0,t]\}$, where $l=r/\alpha+\|a_2+\xi\|$ and $t_\alpha>0$ satisfies $\int_0^{t_\alpha} (a_2(\omega{\cdot}s)+\xi)\, ds\geq \alpha t_\alpha$ for all $\omega\in\Omega$. The existence of such $t_\alpha$ is proved in the proof of Proposition~6.6 of \cite{dno1}. We assume for contradiction that $t_1<\infty$ and define $t_0=\inf\{t<t_1\colon\, w(s)\geq 2r/\alpha$ for all $s\in[t,t_1]\}$. Then, $t_0<t_1-t_\alpha$: otherwise
\begin{equation*}
w(t_1)=w(t_0)+\int_{t_0}^{t_1} \left(-(a_2(\omega_0{\cdot}s)+\xi)+\frac{a_3(\omega_0{\cdot}s)}{w(s)}\right)\, ds<\frac{2r}{\alpha}+lt_\alpha\,,
\end{equation*}
which is not the case. In particular, $w(t)\geq 2r/\alpha$ for $t\in[t_1-t_\alpha,t_1]$, and hence
\begin{equation*}
\begin{split}
w(t_1)&=w(t_1-t_\alpha)-\int_0^{t_\alpha} (a_2(\omega_0{\cdot}(t_1-t_\alpha){\cdot}s)+\xi)\, ds+\int_{t_1-t_\alpha}^{t_1}\frac{a_3(\omega_0{\cdot}s)}{w(s)}\, ds\\ &\leq w(t_1-t_\alpha)-\alpha t_\alpha+\frac{\alpha}{2}t_\alpha<w(t_1-t_\alpha)\,,
\end{split}
\end{equation*}
which contradicts the definition of $t_1$.
\end{proof}
\begin{proposition}\label{prop:4Bpextended} Let $b$ be a continuous primitive of $a_1$. Then,
\begin{enumerate}[label=\rm{(\roman*)}]
\item $\mathrm{sp}(e^ba_2)\subset(0,\infty)$ if and only if \eqref{eq:4fcubic} exhibits the local saddle-node and classical transcritical bifurcations of minimal sets described in Theorem~$\mathrm{\ref{th:clasification}(i)}$, with $\alpha_\lambda$ colliding with $0$ on a residual $\sigma$-invariant set as $\lambda\downarrow\lambda_+$. In particular, this situation holds if $0\not\equiv a_2\geq 0$.
\item $\mathrm{sp}(e^ba_2)\subset(-\infty,0)$ if and only if \eqref{eq:4fcubic} exhibits the local saddle-node and classical transcritical bifurcations of minimal sets described in Theorem~$\mathrm{\ref{th:clasification}(i)}$, with $\beta_\lambda$ colliding with $0$ on a residual $\sigma$-invariant set as $\lambda\downarrow\lambda_+$. In particular, this situation holds if $0\not\equiv a_2\leq 0$.
\item $0\in\mathrm{sp}(e^ba_2)$ if and only if \eqref{eq:4fcubic} exhibits the classical pitchfork bifurcation of minimal sets described in Theorem~$\mathrm{\ref{th:clasification}(ii)}$.
\end{enumerate}
\end{proposition}
\begin{proof} The family of changes of variable $y(t)=e^{-b(\omega{\cdot}t)}x(t)$ takes \eqref{eq:4fcubic} to
\begin{equation}\label{eq:4CP2}
\begin{split}
y'&=-e^{2b(\omega{\cdot}t)}y^3+e^{b(\omega{\cdot}t)}a_2(\omega{\cdot}t)y^2+(a_1(\omega{\cdot}t)-b'(\omega{\cdot}t)+\lambda)y\\
&=-e^{2b(\omega{\cdot}t)}y^3+e^{b(\omega{\cdot}t)}a_2(\omega{\cdot}t)y^2+\lambda y\,.
\end{split}
\end{equation}
The casuistic for \eqref{eq:4fcubic} follows from here, since $\mathcal{N}$ is a minimal set for $\eqref{eq:4CP2}_\lambda$ if and only if $\mathcal{M}=\{(\omega,e^{b(\omega)}x)\colon\,(\omega,x)\in\mathcal{N}\}$ is minimal for $\eqref{eq:4fcubic}_\lambda$. The rest of the proof adapts that of Corollary~6.7 of \cite{dno1}. Theorem~\ref{th:clasification} applied to $f(\omega,y)=-e^{2b(\omega{\cdot}t)}y^3+e^{b(\omega{\cdot}t)}a_2(\omega{\cdot}t)y^2$ shows that $\lambda_-=\lambda_+=0$ is a bifurcation point for \eqref{eq:4CP2}, and that the corresponding bifurcation diagram is that described in its points (i) or (ii). According to Proposition~\ref{prop:newdiagram}, the flow induced by $y'=f(\omega{\cdot}t,y)$, that is, $\eqref{eq:4CP2}_0$, admits just one minimal set if and only if $0\in \mathrm{sp}(e^ba_2)$; two minimal sets, being $\mathcal{M}_0=\Omega\times\{0\}$ the lower one, if and only if $\mathrm{sp}(e^ba_2)\subset(0,\infty)$; and two minimal sets, being $\mathcal{M}_0$ the upper one, if and only if $\mathrm{sp}(e^ba_2)\subset(-\infty,0)$. As said before, this determines the global bifurcation diagram. The last assertions in (i) and (ii) are trivial.
\end{proof}
As a consequence of the previous result, given any strictly positive $a_3$ and any changing-sign $a_2$, we are able to construct $a_1$ with bounded primitive in such a way that \eqref{eq:4fcubic} exhibits the classical pitchfork bifurcation of minimal sets described in Theorem~\ref{th:clasification}(ii).
\begin{proposition}\label{prop:4a2changessign} Let $a_2\in C(\Omega)$ change sign. Then, there exists $a_1\in CP(\Omega)$ such that \eqref{eq:4fcubic} exhibits the classical pitchfork bifurcation described in Theorem~$\mathrm{\ref{th:clasification}(ii)}$.
\end{proposition}
\begin{proof} Let $m\in\mathfrak{M}_\mathrm{erg}(\Omega,\sigma)$ be arbitrarily fixed. We define the nonempty open sets $\Gamma_1=\{\omega\in\Omega\colon\, a_2(\omega)>0\}$ and $\Gamma_2=\{\omega\in\Omega\colon\, a_2(\omega)<0\}$. As $(\Omega,\sigma)$ is minimal, $\mathrm{supp}(m)=\Omega$, so $m(\Gamma_1)>0$ and $m(\Gamma_2)>0$. A suitable application of Urysohn's Lemma provides nonnegative and not identically zero continuous functions $c_1,c_2\colon\Omega\rightarrow\mathbb{R}$ such that $\mathrm{supp}(c_1)\subseteq\Gamma_1$ and $\mathrm{supp}(c_2)\subseteq\Gamma_2$. Then, $\int_\Omega c_1(\omega)a_2(\omega)\, dm>0$ and $\int_\Omega c_2(\omega)a_2(\omega)\, dm<0$, so there exists $\epsilon>0$ such that $\int_\Omega (c_1(\omega)+\epsilon)a_2(\omega)\, dm>0$ and $\int_\Omega (c_2(\omega)+\epsilon)a_2(\omega)\, dm<0$. The density of $C^1(\Omega)$ on $C(\Omega)$ and the strict positiveness of $c_1+\epsilon$ and $c_2+\epsilon$ ensure the existence of strictly positive functions $\hat c_1,\hat c_2\in C^1(\Omega)$ such that $\int_\Omega \hat c_1(\omega)a_2(\omega)\, dm>0$ and $\int_\Omega \hat c_2(\omega)a_2(\omega)\, dm<0$. Therefore, there exists $s\in(0,1)$ such that $\int_\Omega (s\hat c_1(\omega)+(1-s)\hat c_2(\omega))a_2(\omega)\, dm=0$. Since $s\hat c_1(\omega)+(1-s)\hat c_2(\omega)$ is strictly positive, $b(\omega)=\log(s\hat c_1(\omega)+(1-s)\hat c_2(\omega))$ is well defined and $b\in C^1(\Omega)$. So, $\int_\Omega e^{b(\omega)}a_2(\omega)\, dm=0$, and hence $0\in\mathrm{sp}(e^ba_2)$ (see Subsection~\ref{subsec:2measuresandlyapunovexponents}). We take $a_1=b'$ and apply Proposition~\ref{prop:4Bpextended} to complete the proof.
\end{proof}
\subsection{The case of sign-preserving \texorpdfstring{$a_2$}{a2}}\label{subsec:4signpreserving}
In what follows, we will describe some conditions ensuring which one of the bifurcation possibilities of minimal sets described in Theorem~\ref{th:clasification} holds for \eqref{eq:4fcubic}. Our starting point are the functions $a_1,a_3$, with $a_3$ strictly positive. Let $\mathrm{sp}(a_1)=[-\lambda_+,-\lambda_-]$, with $\lambda_-\leq\lambda_+$, be the Sacker and Sell spectrum of $a_1$, let $k_1<k_2$ be such that $k_1\leq a_1(\omega)\leq k_2$ for all $\omega\in\Omega$, and let $0<r_1\leq r_2$ be such that $r_1\leq a_3(\omega)\leq r_2$ for all $\omega\in\Omega$. Since $\mathrm{sp}(a_1)\subseteq [k_1,k_2]$ (see Subsection~\ref{subsec:2measuresandlyapunovexponents}), $\lambda_-+k_1\leq \lambda_++k_1\leq0\leq\lambda_-+k_2\leq\lambda_++k_2$, and the first and last inequalities are strict if $a_1$ has band spectrum. Our goal is to relate the \lq\lq size" of $a_2$ to these six constants in order to get the different diagrams. A remarkable fact is that, in the cases studied this subsection, $a_2$ never changes sign, in contrast with the situation of Proposition \ref{prop:4a2changessign} and those analyzed at the end of Section~\ref{sec:7anotherbifurcation}.

Now, we recall and complete the statement of Proposition~6.5 of \cite{dno1} when applied to our current model \eqref{eq:4fcubic}, which gives a sufficient criterium for the classical pitchfork bifurcation.
\begin{proposition}\label{prop:4classicalcriterium}
\begin{enumerate}[label=\rm{(\roman*)}]
\item {\rm(A criterium ensuring classical pitchfork bifurcation).} If $a_2(\omega)=0$ for all $\omega\in\Omega$, then \eqref{eq:4fcubic} exhibits the classical pitchfork bifurcation of minimal sets described in Theorem~$\mathrm{\ref{th:clasification}(ii)}$.
\item If $a_2(\omega)\geq 0$ (resp. $a_2(\omega)\leq 0$) for all $\omega\in\Omega$, then $\alpha_\lambda$ (resp. $\beta_\lambda$) takes the value $0$ at its continuity points for all $\lambda\leq\lambda_+$.
\end{enumerate}
\end{proposition}
\begin{proof} Proposition~6.5 of \cite{dno1} ensures (i) and (ii) for $\lambda<\lambda_+$. To check (ii) for $\lambda_+$, it suffices to observe that all the alternatives of Theorem~\ref{th:clasification} in which $\alpha_\lambda$ (resp. $\beta_\lambda$) takes the value $0$ at its continuity points for all $\lambda<\lambda_+$, also $\alpha_{\lambda_+}$ (resp. $\beta_{\lambda_+}$) takes the value $0$ at its continuity points.
\end{proof}
The main results of this subsection are stated in Propositions~\ref{prop:4transcriticalcriteriumandnoclassicalpitchfork}, \ref{prop:4notranscritical} and \ref{prop:4generalizedpitchfork}, whose proofs use the next technical results. The first one shows that one of the conditions required in Proposition~\ref{prop:4transcriticalcriteriumandnoclassicalpitchfork} always holds if $a_1$ has band spectrum (in which case $a_1$ is not a constant function).
\begin{lemma}\label{lemma:42espectro} Let $a\in C(\Omega)$. Then, the next three assertions are equivalent: $a$ is nonconstant; $\min_{\omega\in\Omega}a(\w)<\inf \mathrm{sp}(a)$;  $\max_{\omega\in\Omega}a(\w)>\sup\mathrm{sp}(a)$.
\end{lemma}
\begin{proof} Clearly $a$ is nonconstant if $\min_{\omega\in\Omega}a(\w)<\inf \mathrm{sp}(a)$. Now, we assume that $a$ is nonconstant and take
$0<\epsilon<\max_{\omega\in\Omega} a(\omega)-\min_{\omega\in\Omega}a(\omega)$. We define the nonempty open set $U_1=\{\omega\in\Omega\colon\, a(\omega)>\min_{\omega\in\Omega} a(\omega)+\epsilon\}$, and note that $m(U_1)>0$ for all $m\in\mathfrak{M}_\mathrm{erg}(\Omega,\sigma)$, since $\Omega$ is minimal. Let $m_1\in\mathfrak{M}_\mathrm{erg}(\Omega,\sigma)$ satisfy $\inf \mathrm{sp}(a)=\int_\Omega a(\omega)\, dm_1$ (see Subsection~\ref{subsec:2measuresandlyapunovexponents}). Then,
\begin{equation*}
\begin{split}
\inf \mathrm{sp}(a)&=\int_\Omega a(\omega)\, dm_1=\int_{\Omega\setminus U_1} a(\omega)\, dm_1+\int_{U_1} a(\omega)\, dm_1\\
&\geq m_1(\Omega\setminus U_1)\min_{\omega\in\Omega} a(\omega)+m_1(U_1)(\min_{\omega\in\Omega}a(\omega)+\epsilon)>\min_{\omega\in\Omega}a(\omega)\,.
\end{split}
\end{equation*}
This completes the proof of the equivalence of the two first assertions. To check the equivalence of the first and the third ones, we work with $U_2=\{\omega\in\Omega\colon\, a(\omega)<\max_{\omega\in\Omega} a(\omega)-\epsilon\}$.
\end{proof}

\begin{lemma}\label{lemma:4transcriticalcriterium} Let $\lambda+k_1<0$.
\begin{enumerate}[label=\rm{(\roman*)}]
\item If $a_2(\omega)>2\sqrt{r_2(-\lambda-k_1)}$ for all $\omega\in\Omega$, then $\rho_1=\sqrt{(-\lambda-k_1)/r_2}$ is a global strict lower solution of $\eqref{eq:4fcubic}_\lambda$. Consequently, $\tau_\lambda$ admits a minimal set $\mathcal{M}_\lambda^u>\mathcal{M}_0$.
\item If $a_2(\omega)<-2\sqrt{r_2(-\lambda-k_1)}$ for all $\omega\in\Omega$, then $-\rho_1$ is a global strict upper solution of $\eqref{eq:4fcubic}_\lambda$. Consequently, $\tau_\lambda$ admits a minimal set $\mathcal{M}_\lambda^l<\mathcal{M}_0$.
\end{enumerate}
\end{lemma}
\begin{proof} Let us prove (i). We define $g(\rho)=r_2\rho-(\lambda+k_1)/\rho$. Then, $g(\rho_1)=2\sqrt{r_2(-\lambda-k_1)}$, and
\begin{equation*}
\begin{split}
-a_3(\omega)\rho_1^3+a_2(\omega)\rho_1^2+(a_1(\omega)+\lambda)\rho_1&\geq \rho_1^2\left(a_2(\omega)-\left(r_2\rho_1-\frac{\lambda+k_1}{\rho_1}\right)\right)\\
&=\rho_1^2(a_2(\omega)-g(\rho_1))>0
\end{split}
\end{equation*}
for all $\omega\in\Omega$, which proves the first assertion. In turn, this property ensures that $\rho_1<\beta_\lambda$ (see Subsection~\ref{subsec:2compactcoerciveattractor}), which proves the second assertion. The proof of (ii) is analogous.
\end{proof}
\begin{proposition}\label{prop:4transcriticalcriteriumandnoclassicalpitchfork}\begin{enumerate}[label=\rm{(\roman*)}]
\item {\rm(A criterium ensuring saddle-node and transcritical bifurcations).} If $k_1<-\lambda_+$ and
\begin{equation*}
a_2(\omega)>2\sqrt{r_2(-\lambda_--k_1)}\qquad (\text{resp.}\;\;a_2(\omega)<-2\sqrt{r_2(-\lambda_--k_1)}\,)
\end{equation*}
for all $\omega\in\Omega$, then \eqref{eq:4fcubic} exhibits the local saddle-node and transcritical bifurcations of minimal sets described in Theorem~$\mathrm{3.1(i)}$, with $\alpha_\lambda$ (resp. $\beta_\lambda$) colliding with $0$ on a residual $\sigma$-invariant set as $\lambda\downarrow\lambda_+$.
\item {\rm(A criterium precluding classical pitchfork bifurcation).} If $k_1<-\lambda_+$ and
\begin{equation*}
a_2(\omega)>2\sqrt{r_2(-\lambda_+-k_1)}\qquad (resp.\;\;a_2(\omega)<-2\sqrt{r_2(-\lambda_+-k_1)}\,)
\end{equation*}
for all $\omega\in\Omega$, then \eqref{eq:4fcubic} does not exhibit the classical pitchfork bifurcation of minimal sets described in Theorem~$\mathrm{3.1(ii)}$.
\end{enumerate}
\end{proposition}
\begin{proof} (i) Note that if $k_1<-\lambda_+$, then $-\lambda_--k_1\geq-\lambda_+-k_1>0$. There exists $\delta>0$ such that $a_2(\omega)>2\sqrt{r_2(-\lambda_-+\delta-k_1)}$ (resp. $a_2(\omega)<-2\sqrt{r_2(-\lambda_-+\delta-k_1)}$) for all $\omega\in\Omega$. Hence, Lemma~\ref{lemma:4transcriticalcriterium} ensures the existence of a $\tau_{\lambda_--\delta}$-minimal set $\mathcal{M}_{\lambda_--\delta}^u>\mathcal{M}_0$ (resp. $\mathcal{M}_{\lambda_--\delta}^l<\mathcal{M}_0$), and this situation only arises in the stated case of Theorem~3.1(i).

(ii) The existence of a $\tau_{\lambda_+}$-minimal set $\mathcal{M}_{\lambda_+}^u>\mathcal{M}_0$ (resp. $\mathcal{M}_{\lambda_+}^l<\mathcal{M}_0$), ensured by Lemma~\ref{lemma:4transcriticalcriterium}, precludes the situation of Theorem~\ref{th:clasification}(ii).
\end{proof}
Recall that $\lambda_+<k_2$ if $a_1$ has band spectrum: see Lemma~\ref{lemma:42espectro}. The following two results refer to the case that $a_1$ has band spectrum: $\lambda_-<\lambda_+$. (This is ensured in Proposition~\ref{prop:4generalizedpitchfork} by its condition \eqref{eq:4thecondition}).
\begin{proposition}[A criterium ensuring pitchfork bifurcation]\label{prop:4notranscritical} If $\lambda_-<\lambda_+$ and
\begin{equation*}
0\leq a_2(\omega)<\frac{\sqrt{r_1}\,(\lambda_+-\lambda_-)}{\sqrt{\lambda_++k_2}}\qquad \left(\text{resp.}\;\;-\frac{\sqrt{r_1}\,(\lambda_+-\lambda_-)}{\sqrt{\lambda_++k_2}}< a_2(\omega)\leq 0\right)
\end{equation*}
for all $\omega\in\Omega$, then \eqref{eq:4fcubic} does not exhibit the saddle-node and transcritical bifurcations of minimal sets described in Theorem~$\mathrm{\ref{th:clasification}(i)}$.
\end{proposition}
\begin{proof}We reason in the case of positive $a_2$. Let $\hat\beta_\lambda$ be the upper delimiter equilibrium of the global attractor of the skewproduct flow induced by $x'=-a_3(\omega{\cdot}t)x^3+(a_1(\omega{\cdot}t)+\lambda)x$, so that $\hat\beta_\lambda\colon\Omega\rightarrow\mathbb{R}$ is a strictly positive continuous map if $\lambda>\lambda_+$ (see Proposition \ref{prop:4classicalcriterium}(i)). Note that, if $\lambda+k_2\ge 0$ and $\rho>\sqrt{(\lambda+k_2)/r_1}$, then
\begin{equation*}
-a_3(\omega)\rho^3+(a_1(\omega)+\lambda)\rho\leq -r_1\rho^3+(\lambda+k_2)\rho<0\,;
\end{equation*}
that is, $\omega\mapsto\rho$ is a global strict upper solution of $x'=-a_3(\omega{\cdot}t)x^3+(a_1(\omega{\cdot}t)+\lambda)x$. Hence, if $\lambda>\lambda_+$, for any $\omega\in\Omega$, then the $\boldsymbol\upomega$-limit set of $(\omega,\rho)$ contains a minimal set $\mathcal M^u$ which is below of the graph of $\rho$ (see Subsection~\ref{subsec:2equilibria}), and which cannot be $\mathcal M_0$, since $\mathcal M_0$ is repulsive for $\tau_\lambda$. Hence, $\mathcal M^u$ is the graph of $\hat\beta_\lambda$, which ensures that $\hat\beta_\lambda(\omega)\le\rho$ for all $\omega\in\Omega$. It follows easily that $\hat\beta_\lambda(\omega)\le \sqrt{(\lambda+k_2)/r_1}$ for all $\omega\in\Omega$.

We take $\delta>0$ such that, if $\lambda\in[\lambda_+,\lambda_++\delta]$, then $a_2(\omega)\leq \sqrt{r_1}\,(\lambda_+-\lambda_-)/\sqrt{\lambda+k_2}$ for all $\omega\in\Omega$. Then, if $\lambda\in(\lambda_+,\lambda_++\delta]$,
\begin{equation*}
a_2(\omega)\hat\beta_\lambda(\omega)\leq a_2(\omega)\sqrt{(\lambda+k_2)/r_1}\leq \lambda_+-\lambda_-
\end{equation*}
for all $\omega\in\Omega$, and hence
\begin{equation*}
\begin{split}
\hat\beta_\lambda'(\omega)&>-a_3(\omega)\hat\beta_\lambda^3(\omega)+(\lambda_+-\lambda_-)\hat\beta_\lambda(\omega)+(a_1(\omega)+\lambda_-)\hat\beta_\lambda(\omega)\\
&\geq -a_3(\omega)\hat\beta_\lambda^3(\omega)+a_2(\omega)\hat\beta_\lambda^2(\omega)+(a_1(\omega)+\lambda_-)\hat\beta_\lambda(\omega)
\end{split}
\end{equation*}
for all $\omega\in\Omega$: $\hat\beta_\lambda$ is a global strict upper $\tau_{\lambda_-}$-solution, and, in particular, a continuous strong $\tau_{\lambda_-}$-superequilibrium. In addition, Proposition~\ref{prop:4classicalcriterium}(i) ensures that $\hat\beta_{\lambda_+}=\lim_{\lambda\downarrow\lambda_+}\hat\beta_\lambda$ collides with $0$ on a residual set of points $\mathcal{R}_1$.

Let us assume for contradiction that the situation for \eqref{eq:4fcubic} is that described in Theorem~\ref{th:clasification}(i). It follows from Proposition~\ref{prop:4classicalcriterium}(ii) that $\alpha_\lambda$ collides with $0$ on a residual set of points as $\lambda\downarrow\lambda_+$, so Theorem~\ref{th:clasification}(i) ensures that $\beta_{\lambda_-}$ is a continuous strictly positive $\tau_{\lambda_-}$-equilibrium whose graph is the attractive hyperbolic $\tau_{\lambda_-}$-minimal set $\mathcal{M}^u_{\lambda_-}$. Since the Sacker and Sell spectrum of $\mathcal M_0$ for $\tau_{\lambda_-}$ is the nondegenerate interval $[-\lambda_++\lambda_-,0]$, there exists a measure $m$ such that $\int_\Omega (f_x(\omega,0)+\lambda_-)\, dm<0$ (see Subsection~\ref{subsec:2measuresandlyapunovexponents}), and hence Lemma~\ref{lemma:tapaslambdax} ensures that the lower delimiter $\kappa_{\lambda_-}^2$ of the basin of attraction of $\beta_{\lambda_-}$ vanishes on its residual set of continuity points $\mathcal{R}_2$. Let us take $\omega_0\in\mathcal{R}_1\cap\mathcal{R}_2$, a sequence $\{t_n\}\uparrow\infty$ such that $\lim_{n\rightarrow\infty}\omega_0{\cdot}t_n=\omega_0$ and $\lambda\in(\lambda_+,\lambda_++\delta]$ such that $\kappa_{\lambda_-}^2(\omega_0)=0<\hat\beta_\lambda(\omega_0)<\beta_{\lambda_-}(\omega_0)$, which guarantees that $\lim_{n\rightarrow\infty}(u_{\lambda_-}(t_n,\omega_0,\hat\beta_\lambda(\omega_0))-\beta_{\lambda_-}(\omega_0{\cdot}t_n))=0$.  Since $\hat\beta_\lambda$ is a continuous strong $\tau_{\lambda_-}$-superequilibrium, there exist $t_0>0$ and $e_0>0$ such that $u_{\lambda_-}(t,\omega_0,\hat\beta_\lambda(\omega_0))+e\leq\hat\beta_\lambda(\omega_0{\cdot}t)$ for all $t\geq t_0$, and hence $(u_{\lambda_-}(t_n,\omega_0,\hat\beta_\lambda(\omega_0))-\beta_{\lambda_-}(\omega_0{\cdot}t_n))+\beta_{\lambda_-}(\omega_0{\cdot}t_n)+e\leq\hat\beta_\lambda(\omega_0{\cdot}t_n)$ for $n$ large enough. Taking limits as $n\rightarrow\infty$ we obtain $\beta_{\lambda_-}(\omega_0)+e\leq\hat\beta_\lambda(\omega_0)$, a contradiction. This completes the proof.
\end{proof}
\begin{proposition}[A criterium ensuring generalized pitchfork bifurcation]\label{prop:4generalizedpitchfork} If
\begin{equation}\label{eq:4thecondition}
r_1(\lambda_+-\lambda_-)^2+4\,r_2(\lambda_++k_1)(\lambda_++k_2)>0\,,
\end{equation}
and
\begin{equation}\label{eq:4generalizedpitchfork}
\begin{split}
2\sqrt{r_2(-\lambda_+-k_1)}<a_2(\omega)<\frac{(\lambda_+-\lambda_-)\sqrt{r_1}}{\sqrt{\lambda_++k_2}}&\\ \Bigg(\text{resp.}\;\;-\frac{(\lambda_+-\lambda_-)\sqrt{r_1}}{\sqrt{\lambda_++k_2}}< a_2(\omega)<-2\sqrt{r_2(-\lambda_+-k_1)}\Bigg)&
\end{split}
\end{equation}
for all $\omega\in\Omega$, then $\eqref{eq:4fcubic}_\lambda$ exhibits the generalized pitchfork bifurcation of minimal sets described in Theorem~$\mathrm{\ref{th:clasification}(iii)}$, with $\alpha_\lambda$ (resp. $\beta_\lambda$) colliding with $0$ on a residual $\sigma$-invariant set as $\lambda\downarrow\lambda_+$.
\end{proposition}
\begin{proof} Condition \eqref{eq:4thecondition} ensures that $a_1$ has band spectrum and that the intervals in which $a_2$ can take values are nondegenerate. Propositions~\ref{prop:4notranscritical} and \ref{prop:4transcriticalcriteriumandnoclassicalpitchfork}(ii) respectively preclude situations (i) and (ii) of Theorem~\ref{th:clasification}, and Proposition~\ref{prop:4classicalcriterium}(ii) ensures the stated collision property for $\alpha_\lambda$ (resp. for $\beta_\lambda$).
\end{proof}

\subsection{Cases of generalized pitchfork bifurcation}\label{subsec:4cases}
As said after Theorem \ref{th:clasification}, there are autonomous cases presenting either the local saddle-node and classical transcritical bifurcations or the classical pitchfork bifurcation of minimal sets described in cases (i) and (ii) of that theorem. These two possibilities are also the unique ones in nonautonomous examples when $a_1$ has point spectrum, and we have classified them if $a_1\in CP(\Omega)$ in Subsection~\ref{subsec:4CP}, where in addition we have shown simple ways to construct examples fitting in each one of these two situations. In Subsection~\ref{subsec:4signpreserving}, we have observed that cases (i) (with a local nonclassical transcritical bifurcation) and (ii) of Theorem \ref{th:clasification} can also occur when $a_1$ has band spectrum. In addition, Propositions~\ref{prop:4classicalcriterium}(i) and \ref{prop:4transcriticalcriteriumandnoclassicalpitchfork}(i) provide simple ways to construct such examples by choosing a suitable $a_2$ once fixed $a_1$ and $a_3$. In the same line, Proposition \ref{prop:4generalizedpitchfork} establishes conditions ensuring the generalized pitchfork case of Theorem \ref{th:clasification}(iii). But the existence of polynomials satisfying these last conditions is not so obvious.

Therefore, our next objective is to develop systematic ways of constructing third degree polynomials giving rise to families \eqref{eq:4fcubic} for which the global bifurcation diagram is that of Theorem~\ref{th:clasification}(iii). Hence, all the situations described in that theorem actually realize.
\begin{lemma}\label{lemma:7construccioncs}
Let $m_1,\dots,m_n$ be different elements of $\mathfrak{M}_\mathrm{erg}(\Omega,\sigma)$ with $n\geq 1$, and let $0<\epsilon<1$ be fixed. For every $i\in\{1,\dots,n\}$, there exists a continuous $c_i\colon\Omega\rightarrow[0,1]$ with $\min_{\omega\in\Omega}c_i(\omega)=0$ and $\max_{\omega\in\Omega}c_i(\omega)=1$ such that
$c_i·c_j\equiv 0$ and
\begin{equation}
1-\epsilon<\int_\Omega c_i(\omega)\, dm_i\leq1\,,\qquad 0\leq\int_\Omega c_i(\omega)\, dm_j<\epsilon
\end{equation}
for every $i,j\in\{1,\dots,n\}$ with $j\neq i$.
\end{lemma}
\begin{proof}
Let $\epsilon>0$ be fixed. As explained in Remark 1.10 of \cite{johnson1} and Section 6 of Chapter II of \cite{mane1}, there exist disjoint $\sigma$-invariant sets $\Omega_1,\dots,\Omega_n\subseteq\Omega$ such that $m_i(\Omega_j)=\delta_{ij}$ for $i,j\in\{1,\dots,n\}$, where $\delta_{ij}$ is the Kronecker delta: just take $\Omega_i$ as the so-called {\em ergodic component} of $m_i$.
We take compact sets $K_i\subseteq\Omega_i$ such that $m_i(K_i)>1-\epsilon$ for $i\in\{1,\dots,n\}$. Let $d=\min_{1\leq i<j\leq n}d(K_i,K_j)$ and let $U_i$ be an open set such that $K_i\subset U_i\subseteq B(K_i,d/3)$, $m_i(U_i\setminus K_i)<\epsilon$ and $m_j(U_i)=m_j(U_i\setminus K_i)<\epsilon$ for $i,j\in\{1,\dots,n\}$ with $i\neq j$. The choice of $d$ ensures that $U_1,\ldots,U_n$ are pairwise disjoint. Finally, Urysohn's Lemma provides continuous functions $c_i\colon\Omega\rightarrow[0,1]$ with $c_i(\omega)=1$ for all $\omega\in K_i$ and $c_i(\omega)=0$ for all $\omega\notin U_i$, for $i\in\{1,\dots,n\}$. All the statements follow easily.
\end{proof}
\begin{proposition}\label{prop:7linearcombcondition}
Let $m_1,\dots,m_n$ be different elements of $\mathfrak{M}_\mathrm{erg}(\Omega,\sigma)$ with $n\geq2$.
Take $r\ge 1$ and $\epsilon>0$ with
\begin{equation*}
\epsilon<\epsilon_1=\frac{n+2r(n-1)-2\sqrt{r(n-1)\left(r(n-1)+n\right)}}{n^2}\,,
\end{equation*}
(so that $0<\epsilon<1/n$). Let $c_1\ldots,c_n:\Omega\rightarrow\mathbb{R}$ be the functions given by
Lemma~$\mathrm{\ref{lemma:7construccioncs}}$ for $m_1,\ldots,m_n$ and $\epsilon$.
Take constants $\alpha_1\le\alpha_2\leq\dots\le\alpha_n$
with $\alpha_1<0$ and $\alpha_n>0$, and define $a_1=\sum_{i=1}^n \alpha_i c_i$, so that
$\alpha_1\le a_1(\omega)\le\alpha_2$ for all $\omega\in\Omega$.
Then, $a_1$ has band spectrum $\mathrm{sp}(a_1)=[-\lambda_+,-\lambda_-]\subset(\alpha_1,\alpha_2)$ and
\begin{equation}\label{eq:7ecuacionmagica}
(\lambda_+-\lambda_-)^2+4r(\lambda_++\alpha_1)(\lambda_++\alpha_n)>0\,.
\end{equation}
Consequently, if $a_3$ takes values in $[r_1,r_2]$ for $r_1>0$ and $r_2=rr_1$,
and if $a_2$ satisfies \eqref{eq:4generalizedpitchfork} for
$k_1=\alpha_1$ and $k_2=\alpha_n$,
then $\eqref{eq:4fcubic}$ exhibits the generalized pitchfork bifurcation of minimal sets described in Theorem~$\mathrm{\ref{th:clasification}(iii)}$.
\end{proposition}
\begin{proof} It is easy to check that $0<\epsilon_1<1/n$. In addition,
according to Lemma~\ref{lemma:7construccioncs},
\begin{equation*}
\begin{split}
\int_\Omega a_1(\omega)\, dm_1&=\alpha_1\int_\Omega c_1(\omega)\, dm_1+\sum_{i\ge 2}\alpha_i\int_\Omega c_i(\omega)\, dm_1<\alpha_1(1-\epsilon)+(n-1)\alpha_n\epsilon\,,\\
\int_\Omega a_1(\omega)\, dm_n&=\alpha_n\int_\Omega c_n(\omega)\, dm_n+\!\!\!\sum_{i\le n-1}\!\!\!\alpha_i\int_\Omega c_i(\omega)\, dm_n>\alpha_n(1-\epsilon)+(n-1)\alpha_1\epsilon\,.
\end{split}
\end{equation*}
These inequalities ensure that $\alpha_1\leq-\lambda_+<\alpha_1(1-\epsilon)+(n-1)\alpha_n\epsilon$ and $\alpha_n(1-\epsilon)+(n-1)\alpha_1\epsilon<-\lambda_-\leq\alpha_n$.
Then, $\lambda_+-\lambda_->(1-n\epsilon)(\alpha_n-\alpha_1)>0$ (which shows the nondegeneracy of
$\mathrm{sp}(a_1)$), $0>\lambda_++\alpha_1>\epsilon(\alpha_1-(n-1)\alpha_n)$ and $0<\lambda_++\alpha_n\leq\alpha_n-\alpha_1$, which in turn yield
\begin{equation*}
\begin{split}
&(\lambda_+-\lambda_-)^2+4r(\lambda_++\alpha_1)(\lambda_++\alpha_n)\\
&\qquad\qquad>(\alpha_n-\alpha_1)^2(1-n\epsilon)^2+4r\epsilon(\alpha_n-\alpha_1)(\alpha_1-(n-1)\alpha_n)\,.
\end{split}
\end{equation*}
So, it is enough to check that the right-hand side is strictly positive; that is,
\begin{equation*}
\alpha_n((1-n\epsilon)^2-4\epsilon r(n-1))>\alpha_1((1-n\epsilon)^2-4\epsilon r)\,.
\end{equation*}
Since $(1-n\epsilon)^2-4\epsilon r (n-1)\leq (1-n\epsilon)^2-4\epsilon r$ and $\alpha_1<0<\alpha_n$,
it suffices to check that $(1-n\epsilon)^2-4r\epsilon (n-1)>0$. That is, $n^2\epsilon^2-2(n+2r(n-1))\epsilon+1>0$.
And this follows from $\epsilon<\epsilon_1$, since $\epsilon_1$ is the lowest root of the
polynomial.
\par
Since $\alpha_1=\min_{\omega\in\Omega} a_1(\omega)$ and $\alpha_n=\max_{\omega\in\Omega}a_1(\omega)$, Lemma~\ref{lemma:42espectro} ensures that $[-\lambda_+,-\lambda_-]\subset(\alpha_1,\alpha_2)$. Finally, note that all the hypotheses of Proposition~\ref{prop:4generalizedpitchfork} are fulfilled with $k_1=\alpha_1$ and $k_2=\alpha_n$. This proves the last assertion.
\end{proof}

Note that every function $a_1$ constructed by the procedure of Proposition~\ref{prop:7linearcombcondition} takes positive and negative values. But this is not a real restriction to get a generalized pitchfork bifurcation diagram, since that corresponding to $a_1+\mu$ for any constant $\mu\in\mathbb{R}$ is of the same type.

Proposition~\ref{prop:7linearcombcondition} shows that the occurrence of families \eqref{eq:4fcubic} with generalized pitchfork bifurcation diagram only requires the existence of two different ergodic measures. The functions $a_1$ constructed as there indicated are intended to satisfy \eqref{eq:4thecondition}; that is, their extremal Lyapunov exponents are near its maximum and minimum. But in fact this is not a necessary condition for a function $a_1$ to be the first order coefficient of a polynomial giving rise to a generalized pitchfork bifurcation. Theorem~\ref{th:4pitchforkanya1} proves this assertion in the case of a finitely ergodic base flow. Its proof in based on Proposition~\ref{prop:7descompositionCOmega} and Corollary~\ref{corollary:7decomposition}.
\begin{proposition}\label{prop:7descompositionCOmega} Assume that $\mathfrak{M}_\mathrm{erg}(\Omega,\sigma)=\{m_1,\ldots,m_n\}$ with $n\ge 1$.
There exists $\epsilon_2>0$ such that, if $0<\epsilon\leq\epsilon_2$ and $c_1,\ldots,c_n\colon\Omega\rightarrow\mathbb{R}$ are the
functions constructed in Lemma~{\rm\ref{lemma:7construccioncs}} for $m_1,\ldots,m_n$ and $\epsilon$, then
\begin{equation*}
C(\Omega)=\langle c_1,\dots,c_n\rangle\oplus C_0(\Omega)
\end{equation*}
as topological sum of vector spaces, where $C(\Omega)$ is endowed with the uniform topology, given by $\|a\|=\max_{\omega\in\Omega} |a(\omega)|$.
In particular, the Sacker and Sell spectrum of $a\in C(\Omega)$ coincides with that of its projection onto $\langle c_1,\dots,c_n\rangle$.
\end{proposition}
\begin{proof}
Let $M_{n\times n}(\mathbb{R})$ be the linear space of $n\times n$ real matrices, which we endow with the norm $\|C\|_\infty=\max_{1\le i,j\le n}|c_{ij}|$, where $C=\{c_{ij}\}_{1\le i,j\le n}$.
The set of regular $n\times n$ real matrices $GL_n(\mathbb R)$ is an open subset of $M_{n\times n}(\mathbb R)$, and the identity matrix $I$ belongs to $GL_n(\mathbb R)$. Hence, there exists $\epsilon_2\in(0,1)$ such that, if $\|C-I\|_\infty\leq\epsilon_2$, then $C$ is regular. Therefore, if $\epsilon\in(0,\epsilon_2]$, then
the corresponding functions $c_1,\ldots,c_n$ of Lemma \ref{lemma:7construccioncs} provide a regular matrix
\begin{equation*}
C=\left(\begin{array}{ccc}
\int_\Omega c_1(\omega)\, dm_1 & \dots & \int_\Omega c_n(\omega)\, dm_1\\
\vdots & \ddots &\vdots\\
\int_\Omega c_1(\omega)\, dm_n & \dots & \int_\Omega c_n(\omega)\, dm_n\\
\end{array}\right).
\end{equation*}
Let us consider the continuous linear functionals $T_i\colon C(\Omega)\rightarrow\mathbb{R}$, $f\mapsto\int_\Omega f(\omega)\, dm_i$ for $i\in\{1,\dots,n\}$, and note that $\mathrm{Ker}(T_i)$ has codimension 1. Therefore, the codimension of the set $C_0(\Omega)$, which coincides with $\bigcap_{i\in\{1,\dots,n\}}\mathrm{Ker}(T_i)$, is at most $n$. In addition, the linear space $\langle c_1,\dots,c_n\rangle$ has dimension $n$, since the supports of $c_1,\ldots,c_n$ are pairwise disjoint. Let us check that $\langle c_1,\dots,c_n\rangle\cap C_0(\Omega)=\{0\}$: if $c=\sum_{i=1}^n \alpha_i c_i\in C_0(\Omega)$, then $0=\int_\Omega c(\omega)\, dm_j=\sum_{i=1}^n \alpha_i\int_\Omega c_i(\omega)\, dm_j$ for every $j\in\{1,\dots,n\}$. These $n$ equations provide a homogeneous linear system for $\alpha_1,\ldots,\alpha_n$ with regular coefficient matrix $C$; so $\alpha_1=\cdots=\alpha_n=0$ and hence $c=0$. Consequently, $C(\Omega)$ is the algebraic direct sum of $\langle c_1,\dots,c_n\rangle$ and $C_0(\Omega)$. We will check that the projections of $C(\Omega)$ onto each one of the subspaces are continuous, which will complete the proof of the first assertion. Given $a\in C(\Omega)$, its projection $P_{\langle c_1,\dots,c_n\rangle}a=\sum_{i=1}^n \alpha_i c_i$ onto $\langle c_1,\ldots,c_n\rangle$ is given by
\begin{equation*}
\left(\begin{array}{c}
\alpha_1\\
\vdots\\
\alpha_n\end{array}\right)=C^{-1}\left(\begin{array}{c}
\int_\Omega a(\omega)\, dm_1\\
\vdots\\
\int_\Omega a(\omega)\, dm_n\end{array}\right)\,.
\end{equation*}
Therefore, $\|\alpha_i c_i\|=|\alpha_i|\leq n \|C^{-1}\|_{\infty} \|a\|$ for every $i\in\{1,\dots,n\}$, and hence $\|P_{\langle c_1,\dots,c_n\rangle}a\|=\|\sum_{i=1}^n\alpha_i c_i\|\leq n^2 \|C^{-1}\|_\infty\|a\|$. So, $P_{\langle c_1,\dots,c_n\rangle}\colon C(\Omega)\to\langle c_1,\ldots,c_n\rangle$ is continuous. Finally, as $P_{C_0(\Omega)}a=a-P_{\langle c_1,\dots,c_n\rangle}a$, also the projection $P_{C_0(\Omega)}$ is continuous, as asserted. The second assertion is an easy consequence of the first one.
\end{proof}
\begin{corollary}\label{corollary:7decomposition} Assume that $\mathfrak{M}_\mathrm{erg}(\Omega,\sigma)=\{m_1,\ldots,m_n\}$ with $n\ge 2$, and take $r\ge 1$. Let $a_1\in C(\Omega)$ have Sacker and Sell spectrum $sp(a_1)=[-\lambda_+,-\lambda_-]$ with $\lambda_-<0<\lambda_+$. Then, there exist $\tilde a_1\in C(\Omega)$
with Sacker and Sell spectrum $sp(\tilde a_1)=[-\lambda_+,-\lambda_-]$ and $a_1-\tilde a_1\in C_0(\Omega)$, and $k_1,k_2\in\mathbb R$ such that $k_1\le\tilde a_1(\omega)\le k_2$ for all $\omega\in\Omega$ and
\begin{equation}\label{eq:7corollarydecomposition}
(\lambda_+-\lambda_-)^2+4\,r(\lambda_++k_1)(\lambda_++k_2)>0\,.
\end{equation}
\end{corollary}
\begin{proof} We take $\epsilon<\min(\epsilon_1,\epsilon_2)$, with $\epsilon_1$ and $\epsilon_2$ respectively provided by Propositions~\ref{prop:7linearcombcondition} and \ref{prop:7descompositionCOmega}. Let $c_1,\ldots,c_n$ be the functions given by Lemma \ref{lemma:7construccioncs} for $m_1,\ldots,m_n$ and $\epsilon$. Proposition \ref{prop:7descompositionCOmega} provides $\alpha_1,\cdots,\alpha_n$ such that the map $\tilde a_1=\alpha_1 c_1+\cdots+\alpha_nc_n$ satisfies $\mathrm{sp}(\tilde a_1)=\mathrm{sp}(a_1)$.
Hence, since $\lambda_-<0<\lambda_+$, $\tilde a_1$ takes positive and negative values, and therefore there exist $i_1$ with $\alpha_{i_1}<0$ and $i_2$ with $\alpha_{i_2}>0$. The result follows from Proposition \ref{prop:7descompositionCOmega} by reordering the measures $m_1,\ldots,m_n$.
\end{proof}
\begin{theorem}\label{th:4pitchforkanya1} Assume that $\mathfrak{M}_\mathrm{erg}(\Omega,\sigma)=\{m_1,\ldots,m_n\}$ with $n\ge 2$. Let $a_1\in C(\Omega)$ have Sacker and Sell spectrum $\mathrm{sp}(a_1)=[-\lambda_+,-\lambda_-]$ with $\lambda_-<\lambda_+$. Then, there exist strictly positive functions $a_2,a_3\in C(\Omega)$ such that
\eqref{eq:4fcubic} exhibits the generalized pitchfork bifurcation of minimal sets described in Theorem~$\mathrm{\ref{th:clasification}(iii)}$.
\end{theorem}
\begin{proof} We take any strictly positive $\tilde a_3\in C(\Omega)$ and $0<r_1\leq r_2$ with $r_1\leq \tilde a_3(\omega)\leq r_2$ for all $\omega\in\Omega$, and call $r=r_2/r_1$. There is no loss of generality in assuming that $\lambda_-<0<\lambda_+$, since the bifurcation diagrams for $a_1$ and $a_1+\mu$ coincide for any $\mu\in\mathbb{R}$. We associate $\tilde a_1$ to $a_1$ and $r$ by Corollary~\ref{corollary:7decomposition}. Note that there exists $\delta_0>0$ such that, if $|\lambda_+-\mu_+|<\delta_0$ and $|\lambda_--\mu_-|<\delta_0$, then
$\mu_-<\mu_+$ and \eqref{eq:7corollarydecomposition} holds for $\mu_+,\mu_-$ instead of $\lambda_+,\lambda_-$ and $k_1+\delta_0,k_2-\delta_0$ instead of $k_1,k_2$. That is, if $c\in C(\Omega)$ satisfies $|c(\omega)|<\delta_0$ for all $\omega\in\Omega$, then $\text{\rm sp}(\tilde a_1+c)=[-\mu_+,-\mu_-]$ with $\mu_-<\mu_+$, $k_1-\delta_0\le \tilde a_1(\omega)+c(\omega)\le k_2+\delta_0$ for all $\omega\in\Omega$, and
\[
 (\mu_+-\mu_-)^2+4\,r(\mu_++k_1-\delta_0)(\mu_++k_2+\delta_0)>0\,.
\]
Since $CP(\Omega)$ is dense in $C_0(\Omega)$ (see Subsection~\ref{subsec:2boundedprimitive}), there exists $b\in C^1(\Omega)$ such that $\max_{\omega\in\Omega}|a_1(\omega)-\tilde a_1(\omega)-b'(\omega)|<\delta_0$. We take $c=a_1-\tilde a_1-b'$ and call $\text{\rm sp}(\tilde a_1+c)=[-\mu_+,-\mu_-]$. Now, in order to apply Proposition \ref{prop:4generalizedpitchfork}, we take $\tilde a_2\in C(\Omega)$ satisfying $2\sqrt{r_2(\mu_+-k_1-\delta_0)}<\tilde a_2(\omega)<\sqrt{r_1}(\mu_+-\mu_-)/(\sqrt{(\mu_++k_2-\delta_0)})$ for all $\omega\in\Omega$. Hence, the parametric family
\begin{equation}\label{eq:7cambiovariable}
x'=-\tilde a_3(\omega{\cdot}t)x^3+\tilde a_2(\omega{\cdot}t)x^2+(\tilde a_1(\omega{\cdot}t)+c(\omega{\cdot}t)+\lambda)x\,,\quad\omega\in\Omega
\end{equation}
presents a generalized pitchfork bifurcation of minimal sets.
As explained in the proof of Proposition \ref{prop:4Bpextended}, the family of
changes of variables $y(t)=e^{b(\omega{\cdot}t)}x(t)$ takes \eqref{eq:7cambiovariable} to
\begin{equation*}
\begin{split}
y'&=-e^{-2b(\omega{\cdot}t)}\tilde a_3(\omega{\cdot}t)y^3+e^{-b(\omega{\cdot}t)}\tilde a_2(\omega{\cdot}t)y^2+(\tilde a_1(\omega{\cdot}t)+c(\omega{\cdot}t)+b'(\omega{\cdot}t)+\lambda)y\\
&=-e^{-2b(\omega{\cdot}t)}\tilde a_3(\omega{\cdot}t)y^3+e^{-b(\omega{\cdot}t)}\tilde a_2(\omega{\cdot}t)y^2+(a_1(\omega{\cdot}t)+\lambda)y\\
\end{split}
\end{equation*}
without changing the global structure of the bifurcation diagram. That is, the strictly positive functions $a_3=e^{-2b}\tilde a_3$ and $a_2=e^{-b}\tilde a_2$ fulfill the statement.
\end{proof}
\section{Criteria in a more general framework}\label{sec:6generalframework}
The ideas of Subsection~\ref{subsec:4signpreserving} and \ref{subsec:4cases} can be used to construct examples of all the three possible types of global bifurcation diagrams described in Theorem \ref{th:clasification} for families of differential equations of a more general type. Let us consider
\begin{equation}\label{eq:5fgeneral}
 x'=(-a_3(\omega{\cdot}t)+h(\omega{\cdot}t,x))x^3+a_2(\omega{\cdot}t)x^2+(a_1(\omega{\cdot}t)+\lambda)x\,,\quad\omega\in\Omega\,,
\end{equation}
where $a_i\in C(\Omega)$ for $i\in\{1,2,3\}$, $a_3$ is strictly positive, $h\in C^{0,2}(\Omega\times\mathbb{R},\mathbb{R})$, and $h(\omega,0)=0$ for all $\omega\in\Omega$. Throughout the section, we will represent the Sacker and Sell spectrum of $a_1$ as $\text{\rm sp}(a_1)=[-\lambda_+,-\lambda_-]$, with $\lambda_-\le\lambda_+$. In addition, $k_1\le k_2$ represent real constants satisfying $k_1\leq a_1(\omega)\leq k_2$, and $0<r_1\leq r_2$ represent real constants satisfying $r_1\leq a_3(\omega)\leq r_2$ for all $\omega\in\Omega$.

In the line of part of the results of Section \ref{sec:5criteriacubic}, we fix $a_1$, $a_3$ and $h$ satisfying the mentioned hypotheses as well as a fundamental extra condition which relates the behavior of $h$ for small values of $x$ to the properties of $a_1$ and $a_3$, and such that the function
\begin{equation}\label{def:5f}
f(\omega,x)=(-a_3(\omega)+h(\omega,x))x^3+a_2(\omega)x^2+a_1(\omega)x
\end{equation}
is $\mathrm{(Co)}$ and $\mathrm{(SDC)}_*$. The goal is to describe conditions on $a_2$ determining each one of the possible bifurcation cases described in Theorem \ref{th:clasification} for \eqref{eq:5fgeneral}. The function $a_2$ will be sign-preserving under all these conditions. Note the $\mathrm{sp}(a_1)$ is the Sacker and Sell spectrum of $f_x$ on $\mathcal M_0$.
\begin{proposition} \label{prop:5casoshgeneral}
Assume that the function $f$ given by \eqref{def:5f} is $\mathrm{(Co)}$ and $\mathrm{(SDC)}_*$, and that
\begin{itemize}
\item[\rm(H)] there exist $\rho_0>0$ and $0<\epsilon_0<r_1$ such that
$|h(\omega,x)|\leq\epsilon_0$ for all $\omega\in\Omega$ if $|x|\leq\rho_0$,
$\sqrt{(\lambda_++k_2)/(r_1-\epsilon_0)}<\rho_0$, and
$\sqrt{(-\lambda_--k_1)/(r_2+\epsilon_0)}<\rho_0$.
\end{itemize}
We call $s_1=r_1-\epsilon_0$ and $s_0=r_2+\epsilon_0$. Then,
\begin{enumerate}[label=\rm{(\roman*)}]
\item if $a_2(\omega)=0$ for all $\omega\in\Omega$, then \eqref{eq:5fgeneral} exhibits the classical pitchfork bifurcation of minimal sets described in Theorem~$\mathrm{\ref{th:clasification}(ii)}$.
\item If $k_1<-\lambda_+$ and
\begin{equation*}
a_2(\omega)>2\sqrt{s_2(-\lambda_--k_1)}\quad (\text{resp. }a_2(\omega)<-2\sqrt{s_2(-\lambda_--k_1)}\,)
\end{equation*}
for all $\omega\in\Omega$, then \eqref{eq:5fgeneral} exhibits the local saddle-node and transcritical bifurcations of minimal sets described in Theorem~$\mathrm{\ref{th:clasification}(i)}$, with $\alpha_\lambda$ (resp. $\beta_\lambda$) colliding with $0$ on a residual $\sigma$-invariant set as $\lambda\downarrow\lambda_+$.
\item If $k_1<-\lambda_+$ and
\begin{equation*}
a_2(\omega)>2\sqrt{s_2(-\lambda_+-k_1)}\quad(\text{resp. }a_2(\omega)<-2\sqrt{s_2(-\lambda_+-k_1)}\,)
\end{equation*}
for all $\omega\in\Omega$, then \eqref{eq:5fgeneral} does not exhibit the classical pitchfork bifurcation of minimal sets described in Theorem~$\mathrm{\ref{th:clasification}(ii)}$.
\item If $0\leq a_2(\omega)<(\lambda_+-\lambda_-)\sqrt{s_1/(\lambda_++k_2)}$ (resp. $-(\lambda_+-\lambda_-)\sqrt{s_1/(\lambda_++k_2)}<a_2(\omega)\leq0$) for all $\omega\in\Omega$, then \eqref{eq:5fgeneral} does not exhibit the local saddle-node and transcritical bifurcations of minimal sets described in Theorem~$\mathrm{\ref{th:clasification}(i)}$.
\item If
\begin{equation*}
s_1(\lambda_+-\lambda_-)^2+4\,s_2(\lambda_++k_1)(\lambda_++k_2)>0\,,
\end{equation*}
and
\begin{equation*}
\begin{split}
2\sqrt{s_2(-\lambda_+-k_1)}<a_2(\omega)<\frac{(\lambda_+-\lambda_-)\sqrt{s_1}}{\sqrt{\lambda_++k_2}}&\\ \Bigg(\text{resp.}\;\;-\frac{(\lambda_+-\lambda_-)\sqrt{s_1}}{\sqrt{\lambda_++k_2}}< a_2(\omega)<-2\sqrt{s_2(-\lambda_+-k_1)}\Bigg)&
\end{split}
\end{equation*}
for all $\omega\in\Omega$, then \eqref{eq:5fgeneral} exhibits the generalized pitchfork bifurcation of minimal sets described in Theorem~$\mathrm{\ref{th:clasification}(iii)}$.
\end{enumerate}
\end{proposition}
\begin{proof}
(i) This property is proved in Proposition 6.5(ii) of \cite{dno1}.

(ii)-(iii) We define $\rho_\pm=\sqrt{(-\lambda_\pm-k_1)/s_2}$ and observe that $\rho_+\le\rho_-<\rho_0$. Hence, $a_3(\omega)+h(\omega,\rho_\pm)\leq r_2+\epsilon_0=s_2$. This allows us to check the assertions of Lemma~\ref{lemma:4transcriticalcriterium} for $\lambda=\lambda_+$ and $\lambda=\lambda_--\delta$ (for some small $\delta>0$ for which $\sqrt{(-(\lambda_--\delta)-k_1)/s_2}<\rho_0$) and with $r_2$ replaced by $s_2$, just by repeating the same proof. Once this is done, the arguments of Proposition~\ref{prop:4transcriticalcriteriumandnoclassicalpitchfork}(i) and (ii) prove assertions (ii) and (iii).

(iv) The proof is analogous to that of Proposition~\ref{prop:4notranscritical}. We sketch it in the case of positive $a_2$ and refer to Proposition~\ref{prop:4notranscritical} for the details. Condition (H) and the hypothesis in (iv) allow us to take $\delta>0$ such that, if $\lambda\in[\lambda_+,\lambda_++\delta]$, then $\sqrt{(\lambda+k_2)/(r_1-\epsilon_0)}<\rho_0$ and $a_2(\omega)\le(\lambda_+-\lambda_-)\sqrt{s_1/(\lambda+k_2)}$ for all $\omega\in\Omega$. We fix one of these values of $\lambda$ and take $\rho$ satisfying $\sqrt{(\lambda+k_2)/s_1}<\rho\leq\rho_0$ in order to deduce from the bound $a_3(\omega)+h(\omega,\rho)\geq r_1-\epsilon_0=s_1$ (ensured by (H)) that $\rho$ is a strict global upper solution of $x'=(-a_3(\omega{\cdot}t)+h(\omega{\cdot}t,x))x^3+(a_1(\omega{\cdot}t)+\lambda)x$. Let $\hat\beta_\lambda$ be the upper delimiter equilibrium of the global attractor of the corresponding skew-product flow. The next step is checking that $\hat\beta_\lambda(\omega)\leq\sqrt{(\lambda+k_2)/s_1}$ for all $\omega\in\Omega$ if $\lambda\in(\lambda_+,\lambda_++\delta]$. Once this is done, we can deduce that $\hat\beta_\lambda$ is also a strong $\tau_{\lambda_-}$-superequilibrium and that the limit $\hat\beta_{\lambda_+}=\lim_{\lambda\downarrow\lambda_+}\hat\beta_\lambda$ takes the value 0 on a residual set of points. Finally, we assume that \eqref{eq:5fgeneral} is in case (i) of Theorem \ref{th:clasification} and use the previous properties to get a contradiction.

(v) As in Proposition \ref{prop:4generalizedpitchfork}, this follows from (iii), (iv) and Theorem~\ref{th:clasification}.
\end{proof}

We complete this short section by analyzing ways to guarantee the initial hypotheses of Proposition \ref{prop:5casoshgeneral}: $f$ is (Co) and (SDC)$_*$, and (H) holds. The first objectives are achieved, at least, in the following situations:
\begin{proposition} Assume that $(-a_3(\omega)+h(\omega,x))x^3$ is $\mathrm{(Co)}_2$ and $\mathrm{(SDC)}_*$. Then, the function $f$ defined by \eqref{def:5f} is $\mathrm{(Co)}_2$ and $\mathrm{(SDC)}_*$. Moreover, if $(-r+h(\omega,x))x^3$ is $\mathrm{(Co)}$ and $\mathrm{(DC)}$ for some $r<r_1$ (or $r\leq r_1$ if $a_3$ is not constant), then $(-a_3(\omega)+h(\omega,x))x^3$ is $\mathrm{(Co)}_2$ and $\mathrm{(SDC)}_*$.
\end{proposition}
\begin{proof}
Consider $f(\omega,x)=((-a_3(\omega)+h(\omega,x))x^3)+(a_2(\omega)x^2+a_1(\omega)x)$ and recall that any quadratic polynomial is $\mathrm{(DC)}$; that
the sum of a $\mathrm{(SDC)}_*$ function and a $\mathrm{(DC)}$ one is $\mathrm{(SDC)}_*$ (see Section~3 of \cite{dno1}); and that the sum of a $\mathrm{(Co)}_2$ function and a second degree polynomial is $\mathrm{(Co)}_2$.

Let us check the second assertion. Since $(r-a_3(\omega))x^3$ is $\mathrm{(SDC)}_*$ (see Subsection~\ref{subsec2:dconcavity}) and $(-r+h(\omega,x))x^3$ is $\mathrm{(DC)}$, the sum is $\mathrm{(SDC)}_*$. In addition, the $\mathrm{(Co)}$ property of $(-r+h(\omega,x))x^3$ means that $\lim_{|x|\rightarrow\infty} (-r+h(\omega,x))x^2=-\infty$ uniformly on $\Omega$, and hence there exists $\rho>0$ such that $h(\omega,x)<r$ for all $\omega\in\Omega$ and $x\geq\rho$. Consequently,
\begin{equation*}
\limsup_{x\rightarrow\infty}(-a_3(\omega)+h(\omega,x))x\leq\limsup_{x\rightarrow\infty}(-a_3(\omega)+r)x=-\infty
\end{equation*}
for all $\omega\in\Omega$, and hence $\lim_{x\rightarrow\infty}(-a_3(\omega)+h(\omega,x))x=-\infty$ for all $\omega\in\Omega$. The compactness of $\Omega$ ensures that the limit is uniform. Analogously, $\lim_{x\rightarrow-\infty}(-a_3(\omega)+h(\omega,x))x=\infty$ uniformly on $\Omega$. Hence, $(-a_3(\omega)+h(\omega,x))x^3$ is $\mathrm{(Co)}_2$.
\end{proof}

Regarding (H), notice that the included inequality $\sqrt{(-\lambda_--k_1)/(r_2+\epsilon_0)}<\rho_0$ is fulfilled by taking a large enough upper bound $r_2$ for $a_3$ (although the smaller $r_2$ is, the less restrictive the conditions in points (ii), (iii) and (v) of Proposition \ref{prop:5casoshgeneral} are). The following results indicate three ways to get the rest of the conditions in (H). Recall that $h$ is always assumed to belong to $C^{0,2}(\Omega\times\mathbb{R},\mathbb{R})$ and to satisfy $h(\omega,0)=0$ for all $\omega\in\Omega$.
And recall also the meaning of $\lambda_-,\lambda_+,k_1,k_2,r_1$ and $r_2$.

\begin{proposition}
\begin{enumerate}[label=\rm{(\roman*)}]
\item Assume that $a_3$ and $h$ are fixed, take $0<\epsilon_0<r_1$ and $\rho_0>0$ such that $|h(\omega,x)|\leq\epsilon_0$ for all $\omega\in\Omega$ if $|x|\leq\rho_0$, and call $s_1=r_1-\epsilon_0$. Let $k_1<0<k_2$ satisfy $k_2-k_1<\rho_0^2\,s_1$ and $a_1\in C(\Omega)$ satisfy $k_1\leq a_1(\omega)\leq k_2$ for all $\omega\in\Omega$, choose the upper bound $r_2$ for $a_3$ large enough to get $\sqrt{(k_2-k_1)/(r_2+\epsilon_0)}<\rho_0$, and call $s_2=r_2+\epsilon_0$. Then, {\rm(H)} is fulfilled. If, in addition, $(\Omega,\sigma)$ is not uniquely ergodic, then $a_1$ can be chosen to get $\lambda_-<\lambda_+$ and $s_1(\lambda_+-\lambda_-)^2+4\,s_2(\lambda_++k_1)(\lambda_++ k_2)>0$. Consequently, if the function $f$ given by \eqref{def:5f} is $\mathrm{(Co)}$ and $\mathrm{(SDC)}_*$, suitable choices of $a_1$ and $a_2$ provide the three possible global bifurcation diagrams for \eqref{eq:5fgeneral}.
\item Assume that $a_1$ and $h$ are fixed and take $\epsilon_0>0$ and $\rho_0>0$ such that $|h(\omega,x)|\leq\epsilon_0$ for all $\omega\in\Omega$ if $|x|\leq\rho_0$. If $r_1>(\lambda_++k_2)/\rho_0^2+\epsilon_0$ and $a_3\in C(\Omega)$ satisfies $r_1\leq a_3(\omega)$ for all $\omega\in\Omega$, then {\rm(H)} is fulfilled for large enough $r_2$.
\item Assume that $a_1$ and $a_3$ are fixed, and take $\rho_0$ and $m$ with $0<\rho_0<\sqrt{(\lambda_++k_2)/r_1}$ and $0<m< r_1/\rho_0-(\lambda_++k_2)/\rho_0^3$. If $h\in C^{0,2}(\Omega\times\mathbb{R},\mathbb{R})$ satisfies $h(\omega,0)=0$ and $|h_x(\omega,x)|\leq m$ for all $\omega\in\Omega$ if $|x|\leq\rho_0$, then {\rm(H)} is fulfilled with $\epsilon_0=\rho_0 m$ and a large enough $r_2$.
\end{enumerate}
\end{proposition}
\begin{proof} (i) Notice that $\sqrt{(-\lambda_--k_1)/(r_2+\epsilon_0)}\leq\sqrt{(k_2-k_1)/s_2}<\rho_0$ and also $\sqrt{(\lambda_++k_2)/(r_1-\epsilon_0)}\leq\sqrt{(k_2-k_1)/s_1}<\rho_0$. This proves the first assertion.

To prove the second one, we apply Proposition~\ref{prop:7linearcombcondition} with $\alpha_1=k_1$ and $\alpha_n=k_2$. It shows that $a_1$ has band spectrum $[-\lambda_+,-\lambda_-]$ and that $(\lambda_+-\lambda_-)^2+4s(\lambda_++k_1)(\lambda_++k_2)>0$. Lemma~\ref{lemma:42espectro} shows that $\mathrm{sp}(a_1)\subset(k_1,k_2)$. The last assertion in (i) follows from Proposition \ref{prop:5casoshgeneral}.

(ii) If $r_1>(\lambda_++k_2)/\rho_0^2+\epsilon_0$, then $\sqrt{(\lambda_++k_2)/(r_1-\epsilon_0)}<\rho_0$.

(iii) It is clear that $0<\epsilon_0<r_1$, and easy to check that $\sqrt{(\lambda_++k_2)/(r_1-\epsilon_0)}<\rho_0$. A large enough $s_2$ ensures $\sqrt{(-\lambda_--k_1)/(r_2+\epsilon_0)}<\rho_0$.
The equality $h(\omega,x)=\int_0^1 xh_x(\omega,sx)\,ds$ yields $|h(\omega,x)|\leq \rho_0m=\epsilon_0$ for $|x|\leq\rho_0$.
\end{proof}

\section{A second bifurcation problem}\label{sec:7anotherbifurcation}
The ideas and methods developed in \cite{dno1} and in the previous sections of this paper allow us to classify and describe all the possibilities for the bifurcation diagram of a problem different from that analyzed in Sections~\ref{sec:3bifurcationtheorem}, \ref{sec:5criteriacubic} and \ref{sec:6generalframework}, namely
\begin{equation}\label{eq:7anotherbifurcation}
x'=f(\omega{\cdot}t,x)+\mu x^2\,,\quad\omega\in\Omega\,.
\end{equation}
Besides its own interest, this analysis allows us to go deeper in the construction of patterns for the three bifurcation possibilities described in Theorem \ref{th:clasification}, as explained at the end of this section.
\par
As before, $f$ is assumed to be $C^{0,2}(\Omega\times\mathbb{R},\mathbb{R})$ and $\mathrm{(Co)}_2$ (see Subsection \ref{subsec:2compactcoerciveattractor}) with $f(\omega,0)=0$ for all $\omega\in\Omega$, and $\mu\in\mathbb{R}$ is the bifurcation parameter. We expect not to generate risk of confusion while denoting by $\tau_\mu:\mathcal{U}_\mu\subseteq\mathbb{R}\times\Omega\times\mathbb{R}\rightarrow\Omega\times\mathbb{R}$ the local skewproduct flow $(t,\omega,x_0)\mapsto (\omega{\cdot}t,u_\mu(t,\omega,x_0))$ defined from the maximal solutions $\mathcal{I}_{\omega,x_0}^\mu\rightarrow\mathbb{R}$, $t\mapsto u_\mu(t,\omega,x_0)$ of $\eqref{eq:7anotherbifurcation}_\mu$. Note that the set $\mathcal M_0=\Omega\times\{0\}$ is $\tau_\mu$-minimal for all $\mu\in\mathbb{R}$, and that the Sacker and Sell spectrum of $f_x +2\mu x$ on $\mathcal M_0$ (see Subsection~\ref{subsec:2measuresandlyapunovexponents}) is independent of $\mu$. We represent it by $\mathrm{sp}(f_x({\cdot},0))$.

Note that $f+\mu x^2$ is $\mathrm{(Co)}_2$ if $f$ is (unlike what happens with the $\mathrm{(Co)}$ condition).

We will denote by $\mathcal A_\mu$ the corresponding global attractor. As recalled in Subsection~\ref{subsec:2compactcoerciveattractor}, $\mathcal{A}_\mu=\bigcup_{\omega\in\Omega}\big(\{\omega\}\times[\alpha_\mu(\omega),\beta_\mu(\omega)]\big)$, where $\alpha_\mu$ and $\beta_\mu$ are (lower and upper) semicontinuous $\tau_\mu$-equilibria. The following proposition analyzes the behavior of these maps as $\mu$ varies.
\begin{proposition}\label{prop:7properties} Let $f\in C^{0,1}(\Omega\times\mathbb{R},\mathbb{R})$ be $\mathrm{(Co)}_2$. Then,
\begin{enumerate}[label=\rm{(\roman*)}]
\item for any $\omega\in\Omega$, the maps $\mu\mapsto \beta_\mu(\omega)$ and $\mu\mapsto \alpha_\mu(\omega)$ are nondecreasing on $\mathbb{R}$ and they are, respectively, right- and left-continuous. If $\beta_{\mu_0}$ (resp. $\alpha_{\mu_0}$) is strictly positive (resp. strictly negative) for some $\mu_0\in\mathbb{R}$, then $\beta_\mu(\omega)<\beta_{\mu_0}(\omega)$ (resp. $\alpha_\mu(\omega)<\alpha_{\mu_0}(\omega)$) for all $\mu<\mu_0$ and $\omega\in\Omega$, and $\beta_{\mu_0}(\omega)<\beta_\mu(\omega)$ (resp. $\alpha_{\mu_0}(\omega)<\alpha_\mu(\omega)$) for all $\mu>\mu_0$ and $\omega\in\Omega$.
\item $\lim_{\mu\rightarrow-\infty} \alpha_\mu(\omega)=-\infty$ and $\lim_{\mu\rightarrow\infty}\beta_\mu(\omega)=\infty$ uniformly on $\Omega$.
\end{enumerate}
\end{proposition}
\begin{proof}(i) Let $\xi<\mu$. Since $\beta'_\xi(\omega)=f(\omega,\beta_\xi(\omega))+\xi\beta_\xi^2(\omega)\leq f(\omega,\beta_\xi(\omega))+\mu\beta_\xi^2(\omega)$, $\beta_\xi$ is a global lower solution of $\eqref{eq:7anotherbifurcation}_\mu$. Consequently, $\beta_\xi(\omega)\leq\beta_\mu(\omega)$ for all $\omega\in\Omega$ (see Subsection~\ref{subsec:2compactcoerciveattractor}). Moreover, if $\beta_{\mu_0}$ is strictly positive, the previous inequalities are strict. The case of $\alpha_\mu$ follows analogously. To show the one-side continuity, we proceed as in the proof of Theorem~5.5(i) of \cite{dno1}.

(ii) Let $\rho>0$ be fixed. There exists $\mu_\rho>0$ such that $f(\omega,\rho)+\mu_\rho \rho^2\geq 0$ for all $\omega\in\Omega$, and $f(\omega,-\rho)-\mu_\rho (-\rho)^2<0$ for all $\omega\in\Omega$, and hence $\omega\mapsto\rho$ is a global lower $\tau_\mu$-solution and $-\rho$ is a global upper $\tau_{-\mu}$-solution for $\mu\geq\mu_\rho$. This ensures that $\rho\leq\beta_\mu(\omega)$ for all $\omega\in\Omega$ and $-\rho\geq\alpha_{-\mu}(\omega)$ for all $\omega\in\Omega$ if $\mu\ge\mu_\rho$ (see again Subsection~\ref{subsec:2compactcoerciveattractor}), as we wanted to prove.
\end{proof}
The proof of Theorem~\ref{th:xsquareclasification}, which describes the possible bifurcation diagrams for \eqref{eq:7anotherbifurcation}, requires the next technical result, similar to Proposition~4.4 of \cite{dno1}.

\begin{proposition}\label{prop:7sumofexponents} Let $f\in C^{0,2}(\Omega\times\mathbb{R},\mathbb{R})$ be $\mathrm{(DC)}$, let us fix $m\in\mathfrak{M}_\mathrm{erg}(\Omega,\sigma)$, $\lambda_0>0$ and $\nu<\mu$ (resp. $\mu<\nu$), let $\kappa_\mu^1\colon\Omega\rightarrow\mathbb{R}$ be an $m$-measurable $\tau_\mu$-equilibrium, and let $\kappa_\nu^2\colon\Omega\rightarrow\mathbb{R}$ be an $m$-measurable equilibrium of $x'=f(\omega{\cdot}t,x)-\lambda_0x+\nu x^2$ such that $0<\kappa_\nu^2(\omega)<\kappa_\mu^1(\omega)$ (resp. $\kappa_\mu^1(\omega)<\kappa_\nu^2(\omega)<0$) for $m$-a.e. $\omega\in\Omega$. Then,
\begin{equation*}
\int_\Omega \big(f_x(\omega,\kappa_\mu^1(\omega))+2\mu\kappa_\mu^1(\omega)\big)\, dm+\int_\Omega\big(f_x(\omega,\kappa_\nu^2(\omega))-\lambda_0+2\nu\kappa_\nu^2(\omega)\big)\, dm<0\,.
\end{equation*}
\end{proposition}
\begin{proof} We define the $m$-a.e. positive (resp. negative) function $k(\omega)=\kappa_\mu^1(\omega)-\kappa_\nu^2(\omega)$ and $F(\omega,y)=\int_0^1 f_x(\omega,sy+\kappa_\nu^2(\omega))\, ds$. Then, it is not hard to check that
\begin{equation*}
\frac{k'(\omega{\cdot}t)}{k(\omega{\cdot}t)}=F(\omega{\cdot}t,k(\omega{\cdot}t))+(\mu-\nu)\frac{(\kappa_\mu^1(\omega{\cdot}t))^2}{k(\omega{\cdot}t)}+\nu(\kappa_\mu^1(\omega{\cdot}t)+\kappa_\nu^2(\omega{\cdot}t))+\lambda_0\frac{\kappa_\nu^2(\omega{\cdot}t)}{k(\omega{\cdot}t)}\,.
\end{equation*}
Notice that $\kappa_\mu^1$, $\kappa_\nu^2$ and $k$ are bounded (see Subsection~\ref{subsec:2compactcoerciveattractor}), that $k'(\omega)/k(\omega)-F(\omega,k(\omega))-\nu(\kappa_\mu^1(\omega)+\kappa_\nu^2(\omega))$ is strictly positive for $m$-a.e. $\omega\in\Omega$, and that $F(\omega,k(\omega))$ and $\nu(\kappa_\mu^1(\omega)+\kappa_\nu^2(\omega))$ are in $L^1(\Omega)$. Then, Birkoff's Ergodic Theorem (see Theorem~1 in Section 2 of Chapter 1 of \cite{sinai1} and Proposition~1.4 of \cite{johnson1}) yields
\begin{equation}\label{eq:7betprop}
\begin{split}
&0=\int_\Omega F(\omega,k(\omega))\, dm+(\mu-\nu)\int_\Omega\frac{(\kappa_\mu^1(\omega))^2}{k(\omega)}\, dm\\
&\qquad\qquad+\nu\int_\Omega(\kappa_\mu^1(\omega)+\kappa_\nu^2(\omega))\, dm+\lambda_0\int_\Omega\frac{\kappa_\nu^2(\omega)}{k(\omega)}\, dm\,.
\end{split}
\end{equation}
The d-concavity of $f$ and Taylor's Theorem ensure that $F(\omega,k(\omega))-F(\omega,0)\geq k(\omega)F_y(\omega,k(\omega))$ for all $\omega\in\Omega$. Now, we derive $yF(\omega,y)=f(\omega,y+\kappa_\nu^2(\omega))-f(\omega,\kappa_\nu^2(\omega))$ with respect to $y$ and evaluate at $y=k(\omega)$, obtaining $F(\omega,k(\omega))+k(\omega)F_y(\omega,k(\omega))=f_x(\omega,\kappa_\mu^1(\omega))$. Since $F(\omega,0)=f_x(\omega,\kappa_\nu^2(\omega))$, the last equality and the previous inequality yield $f_x(\omega,\kappa_\mu^1(\omega))+f_x(\omega,\kappa_\nu^2(\omega))\leq 2F(\omega,k(\omega))$. In turn, this inequality and \eqref{eq:7betprop} yield
\begin{equation*}
\begin{split}
&\int_\Omega (f_x(\omega,\kappa_\mu^1(\omega))+2\mu\kappa_\mu^1(\omega))\, dm+\int_\Omega(f_x(\omega,\kappa_\nu^2(\omega))-\lambda_0+2\nu \kappa_\nu^2(\omega))\, dm\\
&\qquad\leq 2\int_\Omega F(\omega,k(\omega))\, dm+2\int_ \Omega(\mu\kappa_\mu^1(\omega)+\nu\kappa_\nu^2(\omega))\,dm-\lambda_0\\
&\qquad=-2(\mu-\nu)\int_ \Omega\frac{\kappa_\mu^1(\omega)\kappa_\nu^2(\omega)}{k(\omega)}\, dm-\lambda_0\int_\Omega\frac{\kappa_\mu^1(\omega)+\kappa_\nu^2(\omega)}{k(\omega)}\, dm<0\,,
\end{split}
\end{equation*}
which proves the statement in any of the two stated situations.
\end{proof}
\begin{theorem}\label{th:xsquareclasification} Let $f\in C^{0,2}(\Omega\times\mathbb{R},\mathbb{R})$ be $\mathrm{(Co)}_2$ and $\mathrm{(SDC)}_*$.
\begin{enumerate}[label=\rm{(\roman*)}]
\item {\rm(No bifurcation).} If $\mathrm{sp}(f_x(\cdot,0))\subset(0,\infty)$, then for all $\mu\in\mathbb{R}$ there exist three hyperbolic $\tau_\mu$-minimal sets $\mathcal{M}_\mu^l<\mathcal{M}_0=\Omega\times\{0\}<\mathcal{M}_\mu^u$, where $\mathcal{M}_\mu^l$ and $\mathcal{M}_\mu^u$ are respectively attractive hyperbolic copies of the base and given by the graphs of $\alpha_\mu$ and $\beta_\mu$, $\mathcal{M}_0$ is a repulsive hyperbolic copy of the base, and $\lim_{\mu\rightarrow\infty}\alpha_\mu(\omega)=\lim_{\mu\rightarrow-\infty}\beta_\mu(\omega)=0$ uniformly on $\Omega$.
\item {\rm(Two local saddle-node bifurcations).} If $\mathrm{sp}(f_x(\cdot,0))\subset(-\infty,0)$, then $\mathcal{M}_0$ is an attractive hyperbolic $\tau_\mu$-copy of the base for all $\mu\in\mathbb{R}$ and there exist $\mu_1<\mu_2$ such that: for all $\mu<\mu_1$ (resp. $\mu>\mu_2$), there exist three hyperbolic $\tau_\mu$-copies of the base $\mathcal{M}_\mu^l<\mathcal{N}_\mu<\mathcal{M}_0$ (resp. $\mathcal{M}_0<\mathcal{N}_\mu<\mathcal{M}_\mu^u$) which are hyperbolic copies of the base, given by the graphs of $\alpha_\mu<\kappa_\mu<0$ (resp. $0<\kappa_\mu<\beta_\mu$), where $\mathcal{M}_\mu^l$ (resp. $\mathcal{M}_\mu^u$) is hyperbolic attractive and $\mathcal{N}_\mu$ is hyperbolic repulsive, and $\mu\mapsto \kappa_\mu$ is strictly decreasing on $(-\infty,\mu_1)$ (resp. $(\mu_2,\infty)$); the graphs of $\alpha_\mu$ and $\kappa_\mu$ (resp. $\kappa_\mu$ and $\beta_\mu$) collide on a residual $\sigma$-invariant set as $\mu\uparrow\mu_1$ (resp. $\mu\downarrow\mu_2$), giving rise to a nonhyperbolic minimal set $\mathcal{M}^l_{\mu_1}$ (resp. $\mathcal{M}^u_{\mu_2}$); for $\mu\in(\mu_1,\mu_2)$, $\mathcal{A}_\mu=\mathcal{M}_0$; and $\lim_{\mu\rightarrow\pm\infty}\kappa_\mu(\omega)=0$ uniformly on $\Omega$.
\item {\rm(Weak generalized transcritical bifurcation).} If $0\in \mathrm{sp}(f_x(\cdot,0))$, then $\mathcal{M}_0$ is a nonhyperbolic $\tau_\mu$-copy of the base for all $\mu\in\mathbb{R}$. In addition, there exist $\mu_1\leq\mu_2$ such that: for all $\mu<\mu_1$ (resp. $\mu>\mu_2$) there exist exactly two $\tau_\mu$-minimal sets $\mathcal{M}_\mu^l<\mathcal{M}_0$ (resp. $\mathcal{M}_0<\mathcal{M}_\mu^u$), where $\mathcal{M}_\mu^l$ (resp. $\mathcal{M}_\mu^u$) is an attractive hyperbolic copy of the base given by the graph of $\alpha_\mu$ (resp. $\beta_\mu$); and if $\mathcal{M}_0$ is the unique $\tau_{\mu_1}$-minimal set (resp. $\tau_{\mu_2}$-minimal set), then $\alpha_\mu$ and $0$ (resp. $0$ and $\beta_\mu$) collide on a residual $\sigma$-invariant set as $\mu\uparrow\mu_1$ (resp. $\mu\downarrow\mu_2$). In addition, if $0\neq\inf\mathrm{sp}(f_x(\cdot,0))$, then $\mu_1<\mu_2$ and $\mathcal{M}_0$ is the unique $\tau_\mu$-minimal set for any $\mu\in(\mu_1,\mu_2)$; and if $0=\inf\mathrm{sp}(f_x(\cdot,0))$, then $\mathcal{M}_0$ is the unique $\tau_{\mu}$-minimal set for all $\mu\in[\mu_1,\mu_2]$.
\end{enumerate}
\end{theorem}
\begin{proof}
We call $\mathcal M_\mu^l$ and $\mathcal M_\mu^u$ the lower and upper $\tau_\mu$-minimal sets, defined as in Section~\ref{sec:3bifurcationtheorem}, and recall that they are attractive if they are hyperbolic, in which case they respectively coincide with the graphs of the continuous maps $\alpha_\mu$ and $\beta_\mu$.

(i) Since every Lyapunov exponent of $\mathcal{M}_0$ is strictly positive for any $\mu\in\mathbb{R}$ (see Subsection~\ref{subsec:2measuresandlyapunovexponents}), $\mathcal{M}_0$ is a repulsive hyperbolic $\tau_\mu$-minimal set for every $\mu\in\mathbb{R}$. Consequently, there exist three different hyperbolic $\tau_\mu$-minimal sets $\mathcal{M}_\mu^l<\mathcal{M}_0<\mathcal{M}_\mu^u$ (see Subsection~\ref{subsec2:dconcavity}), with $\mathcal{M}_\mu^l$ and $\mathcal{M}_\mu^u$ given respectively by the graphs of $\alpha_\mu<\beta_\mu$, which are continuous. The hyperbolic continuation of minimal sets (see Theorem~3.8 of \cite{potz}) guarantees the continuity of the maps $\mathbb{R}\rightarrow C(\Omega)$, $\mu\mapsto\beta_\mu$ and $\mathbb{R}\rightarrow C(\Omega)$, $\mu\mapsto\alpha_\mu$ in the uniform topology.

To check that $\lim_{\mu\rightarrow-\infty}\beta_\mu(\omega)=0$ uniformly on $\Omega$, we take any $\epsilon>0$, use coercivity to find $\rho>\epsilon>0$ such that $f(\omega,x)\leq 0$ for all $x\geq\rho$ and $\omega\in\Omega$, and choose $\mu_\epsilon<0$ such that $f(\omega,x)+\mu_\epsilon x^2\leq 0$ for all $x\in[\epsilon,\rho]$ and $\omega\in\Omega$. Then, $f(\omega,x)+\mu_\epsilon x^2\leq 0$ for all $x\geq\epsilon$ and $\omega\in\Omega$. According to Theorem~5.1(i) of \cite{dno1}, $\beta_\mu(\omega)\leq\epsilon$ for all $\omega\in\Omega$ if $\mu\leq\mu_\epsilon$, which proves the assertion. The argument is analogous for $\alpha_\mu$.

(ii) Since every Lyapunov exponent of $\mathcal{M}_0$ is strictly negative, $\mathcal M_0$ is always an attractive hyperbolic $\tau_\mu$-minimal set for all $\mu\in\mathbb R$. We fix $\rho>0$ and look for $\mu^+_\rho$ such that $f(\w,\rho)+\mu^+_\rho\rho^2>0$: if $\mu\ge\mu^+_\rho$, then $\rho$ is a global strict lower solution of \eqref{eq:7anotherbifurcation}$_\mu$, and hence its graph is strictly below (resp.~strictly above) a minimal set $\mathcal M^u_\mu$ (resp.~$\mathcal N_\mu$) contained in the $\boldsymbol\omega$-limit (resp.~$\boldsymbol\alpha$-limit) of any point $(\w,\rho)$ (see Subsection \ref{subsec:2equilibria}). The attractive hyperbolicity of $\mathcal{M}_0$ precludes $\mathcal{N}_\mu=\mathcal{M}_0$. That is, there exist three distinct $\tau_\mu$-minimal sets, so they are hyperbolic copies of the base, with $\mathcal M^u_\mu$ attractive (and given by the graph of $\beta_\mu$) and $\mathcal N_\mu$ repulsive. A similar argument works for $-\rho$ and $\mu\le\mu^-_\rho$ if $f(\w,-\rho)+\mu^-_\rho\rho^2<0$, providing three different hyperbolic $\tau_\mu$-copies of the base $\mathcal M^l_\mu<\mathcal N_\mu<\mathcal M_0$, with $\mathcal M^l_\mu$ attractive (and given by the graph of $\alpha_\mu$) and $\mathcal N_\mu$ repulsive. Let $\kappa_\mu$ be $\tau_\mu$-equilibrium whose graph is $\mathcal N_\mu$, both for $\mu\ge\mu_\rho^+$ and $\mu\le\mu_\rho^-$. Since the initially fixed $\rho>0$ is as small as desired, $\lim_{\mu\rightarrow\pm\infty}\kappa_\mu(\omega)=0$ uniformly on $\Omega$.

Let us define $\mathcal I_1=\{\mu\in\mathbb{R}\colon\,$ $\mathcal M_\xi^l$ is hyperbolic and $\mathcal M_\xi^l<\mathcal M_0$ for all $\xi<\mu\}$ and $\mathcal I_2=\{\mu\in\mathbb{R}\colon\,$ $\mathcal M_\xi^u$ is hyperbolic and $\mathcal M_\xi^u>\mathcal M_0$ for all $\xi>\mu\}$, and observe that $\mu_\rho^-\in\mathcal I_1$ and $\mu_\rho^+\in\mathcal I_2$. We also define $\mu_1=\sup\mathcal I_1$ and $\mu_2=\inf\mathcal I_2$ and note that $\mu_1\notin\mathcal I_1$ and $\mu_2\notin\mathcal I_2$: $\mathcal I_1=(-\infty,\mu_1)$ and $\mathcal I_2=(\mu_2,\infty)$.  The existence of at most three $\tau_\mu$-minimal sets for any $\mu$ ensures $\mu_1\le\mu_2$, and hence they are finite. If $\mu\in\mathcal I_1$ (resp.~$\mu\in\mathcal I_2$) then $\mathcal M^l_\mu$ (resp.~$\mathcal M^u_\mu$) is attractive and coincides with the graph of $\alpha_\mu$ (resp. $\beta_\mu$), which is continuous, and there exists a repulsive hyperbolic copy of the base $\mathcal N_\mu$ for $\tau_\mu$ with $\mathcal M_\xi^l<\mathcal N_\mu<\mathcal M_0$ (resp. $\mathcal M_0<\mathcal N_\mu<\mathcal M_\xi^l$): see Subsection \ref{subsec2:dconcavity}. Let $\kappa_\mu$ be the continuous $\tau_\mu$-equilibrium whose graph is $\mathcal N_\mu$ for $\mu\in\mathcal I_1\cup\mathcal I_2$. As in (i), the hyperbolic continuation of minimal sets guarantees the continuity with respect to the uniform topology of the maps $\mu\mapsto\alpha_\mu,\kappa_\mu$ on $\mathcal I_1$ and $\mu\mapsto\kappa_\mu,\beta_\mu$ on $\mathcal I_2$. In addition, if $\xi_1<\xi_2<\mu_1$ or $\mu_2<\xi_1<\xi_2$, then $\kappa_{\xi_1}$ is a strong time-reversed $\tau_{\xi_2}$-superequilibrium, and hence the $\boldsymbol{\upalpha}$-limit set for $\tau_{\xi_2}$ of a point $(\w,\kappa_{\xi_1}(\w))$ contains a $\tau_{\xi_2}$-minimal set which is strictly below the graph $\mathcal N_{\xi_1}$ of $\kappa_{\xi_1}$, and which (see again Subsection~\ref{subsec:2equilibria}), as above, must coincide with $\mathcal N_{\xi_2}$. That is, $\mu\mapsto\kappa_\mu(\omega)$ is strictly decreasing on $\mathcal I_1$ and $\mathcal I_2$ for all $\omega\in\Omega$.

The maps $\alpha_{\mu_1}(\omega)=\lim_{\mu\uparrow\mu_1}\alpha_\mu(\omega)$ and $\kappa_{\mu_1}(\omega)=\lim_{\mu\uparrow\mu_1}\kappa_\mu(\omega)$ satisfy $\alpha_{\mu_1}\le \kappa_{\mu_1}<0$ and respectively define lower and upper semicontinuous $\tau_{\mu_1}$-equilibria which must coincide on the residual set of its common continuity points: otherwise they would define two different $\tau_{\mu_1}$-minimal sets $\mathcal M^l_{\mu_1}<\mathcal N_{\mu_1}<\mathcal M_0$ (see Proposition 2.3 of \cite{dno1}), so $\mathcal M^l_{\mu_1}$ would be hyperbolic attractive and $\mu_1\in\mathcal I_1$, which is not the case.
Therefore, $\mu_1$ is a local saddle-node bifurcation point of minimal sets and a point of discontinuity of the global attractor. In addition, the unique $\tau_{\mu_1}$-minimal set $\mathcal M_{\mu_1}^l$ that they define is nonhyperbolic, which also follows from $\mu_1\notin\mathcal I_1$. Similar arguments apply to $\mu_2$.
In addition, notice that $\mu_1<\mu_2$, as otherwise there would exist three $\tau_{\mu_1}$-minimal sets, contradicting the nonhyperbolicity of $\mathcal{M}_{\mu_1}^l$.
\par

Let $\mu\in(\mu_1,\mu_2)$, and let $\w_0$ be a common continuity point of $\alpha_\mu$, $\alpha_{\mu_1}$ and $\kappa_{\mu_1}$, so that $\alpha_{\mu_1}(\omega_0)=\kappa_{\mu_1}(\omega_0)$. Since $\kappa_{\mu_1}$ is a strong $\tau_\mu$-subequilibrium and $\alpha_{\mu_1}\le\alpha_\mu$, there exists $t_0>0$ such that
$\kappa_{\mu_1}(\omega_0{\cdot}·t_0)<u_\mu(t_0,\omega_0,\alpha_{\mu_1}(\omega_0))\leq\alpha_\mu(\omega_0{\cdot}t_0)$. Let us call $\w_1=\w_0{\cdot}t_0$ and assume for contradiction that $\alpha_\mu(\omega_0)<0$ or, equivalently, $\alpha_\mu(\omega_1)<0$. The previous property combined with $\lim_{\xi\to-\infty}\kappa_\xi(\omega_1)=0$ and the continuity of $\xi\mapsto\kappa_\xi(\omega_1)$ on $\mathcal I_1$ provides $\alpha_\mu(\omega_1)=\kappa_{\xi_1}(\omega_1)$ for a $\xi_1\in\mathcal I_1$. Since $\kappa_{\xi_1}$ is a strong (and continuous) $\tau_\mu$-subequilibrium, there exists $s_1>0$ and $e_1>0$ such that $\alpha_\mu(\omega_1{\cdot}t)>\kappa_{\xi_1}(\omega_1{\cdot}t)+e_1$ for all $t\ge s_1$ (see Subsection \ref{subsec:2equilibria}), and we obtain the contradiction $\alpha_\mu(\omega_1)\ge \kappa_{\xi_1}(\omega_1)+e_1$ by taking $\{t_n\}\uparrow\infty$ with $\omega_1=\lim_{n\to\infty}\omega_1{\cdot}t_n$. This means that $\mathcal M^l_\mu=\mathcal M_0$, and a similar argument shows $\mathcal M^u_\mu=\mathcal M_0$. Consequently, $\mathcal M_0$ is the unique $\tau_\mu$ minimal set for $\mu\in(\mu_1,\mu_2)$ and, since it is hyperbolic attractive, Theorem~3.4 of \cite{dno1} proves that $\mathcal A_\mu=\mathcal M_0$.

(iii) In this case, $\mathcal M_0$ is nonhyperbolic for all $\mu\in\mathbb{R}$ (see Subsection~\ref{subsec:2measuresandlyapunovexponents}), and this fact precludes the existence of three $\tau_\mu$-minimal sets (see Subsection \ref{subsec2:dconcavity}). Let us check that $\mathcal M^u_\mu$ is an attractive hyperbolic copy of the base if $\mu$ is large enough. To this end, we choose $\lambda_0>0$ such that the Sacker and Sell spectrum of $f_x(\omega,0)-\lambda_0$ is contained in $(-\infty,0)$. Consequently, the bifurcation diagram of minimal sets of
\begin{equation}\label{eq:7auxiliar}
x'=f(\omega{\cdot}t,x)-\lambda_0x+\nu x^2
\end{equation}
with respect to $\nu$ is that described in (ii). Let $\nu_2$ be the upper bifurcation point, fix $\nu>\nu_2$, and let $\hat\kappa_\nu$ and $\hat\beta_\nu$ be the equilibria giving rise to the repulsive and attractive hyperbolic copies of the base for $\eqref{eq:7auxiliar}_\nu$. Then, if $\mu>\nu$,
\begin{equation*}
f(\omega,\hat\beta_\nu(\omega))+\mu\hat\beta_\nu^2(\omega)>f(\omega,\hat\beta_\nu(\omega))-\lambda_0\hat\beta_\nu(\omega)+\nu\hat\beta_\nu^2(\omega)\,,
\end{equation*}
 which means that $\hat\beta_\nu$ is a global strict lower $\tau_\mu$-solution and hence that $\hat\kappa_\nu(\omega)<\hat\beta_\nu(\omega)<\beta_\mu(\omega)$ for all $\omega\in\Omega$ (see Subsection~\ref{subsec:2compactcoerciveattractor}). The definition of $\mathcal M^u_\mu$ ensures that it is above the graph of $\hat\beta_\nu$, and hence $\hat\kappa_\nu(\omega)<\gamma_\mu(\omega)$ for any $\tau_\mu$-equilibrium $\gamma_\mu\colon\Omega\rightarrow\mathbb{R}$ with graph contained in $\mathcal M^u_\mu$. Since $\hat\kappa_\nu$ defines a repulsive hyperbolic copy of the base of $\eqref{eq:7auxiliar}_\nu$, Proposition~\ref{prop:7sumofexponents} ensures that $\int_\Omega (f_x(\omega,\gamma_\mu(\omega))+2\mu\gamma_\mu(\omega))\, dm<0$ for all $m\in\mathfrak{M}_\mathrm{erg}(\Omega,\sigma)$. This means that $\mathcal{M}_\mu^u$ is an attractive hyperbolic minimal set for all $\mu>\nu_2$. An analogous argument shows that $\mathcal{M}_\mu^l$ is an attractive hyperbolic minimal set if $-\mu$ is large enough.

Now, we define $\mathcal I_1$, $\mathcal I_2$, $\mu_1$ and $\mu_2$ as in (ii). Note that $\mu_1\le\mu_2$: otherwise there would be three $\tau_{\mu_1}$-minimal sets. The definition of these two values of the parameter shows the first assertion in (iii), concerning $\mu\notin[\mu_1,\mu_2]$. In addition, the definition on $\mathcal M^l_\mu$ (resp.~$\mathcal M^u_\mu$) shows that, if $\mathcal{M}_0$ is the unique $\tau_{\mu_1}$-minimal set (resp. $\tau_{\mu_2}$-minimal set), then $\alpha_\mu$ and $0$ (resp. $0$ and $\beta_\mu$) collide on a residual $\sigma$-invariant set as $\mu\uparrow\mu_1$ (resp. $\mu\downarrow\mu_2$).

Let us prove that $\mathcal M_0$ is the unique $\tau_\mu$-minimal set if $\mu\in[\mu_1,\mu_2]$ in the case $\inf\mathrm{sp}(f_x(\cdot,0))=0$, which means that $\int_\Omega f_x(\w,0)\,dm\ge 0$ for all $m\in \mathfrak{M}_\mathrm{erg}(\Omega,\sigma)$. In this case, any $\tau_\mu$-minimal set different from $\mathcal{M}_0$ is hyperbolic attractive for any $\mu\in\mathbb R$: Proposition~4.4 of \cite{dno1} ensures that all the Lyapunov exponents of such a minimal set are strictly negative, from where the assertion follows. Hence, $\mathcal{M}_{\mu_1}^l$ and $\mathcal{M}_{\mu_2}^u$ coincide with $\mathcal{M}_0$ to not contradict the definition of $\mu_1$ and $\mu_2$. The monotonicity properties of $\alpha_\mu$ and $\beta_\mu$ established in Proposition \ref{prop:7properties} ensure that $\mathcal M^l_\mu=\mathcal M_0$ for all $\mu\ge\mu_1$ and $\mathcal M^u_\mu=\mathcal M_0$ for all $\mu\le\mu_1$, and this proves the assertion.

In the rest of the proof, we work in the remaining case, $\inf\mathrm{sp}(f_x(\cdot,0))<0$, which implies the existence of a measure $m\in \mathfrak{M}_\mathrm{erg}(\Omega,\sigma)$ with $\int_\Omega f_x(\w,0)\,dm<0$.

Let us prove that $\mu_1<\mu_2$. We fix $m\in\mathfrak{M}_\mathrm{erg}(\Omega,\sigma)$ with $\int_\Omega f_x(\w,0)\,dm<0$. For $\mu>\mu_2$, Lemma \ref{lemma:tapaslambdax} ensures that $\int_\Omega (f_x(\omega,\kappa^2_\mu(\omega))+2\mu\kappa^2_\mu(\omega))\,dm>0$, where $\kappa^2_\mu$ is the $m$-measurable $\tau_\mu$-equilibrium providing the lower delimiter of the basin of attraction of the graph of $\beta_\mu$, which satisfies $0\le \kappa_\mu^2<\beta_\mu$. Reasoning as in Theorem~\ref{th:clasification}, we check that $\mu\mapsto\kappa^2_\mu(\omega)$ is nonincreasing on $\mathcal I_2$ for all $\omega\in\Omega$. Hence, there exists the limit $\kappa^2_{\mu_2}=\lim_{\mu\downarrow\mu_2}\kappa^2_\mu\ge 0$. Lebesgue's Convergence Theorem ensures that $\int_\Omega(f_x(\omega,\kappa^2_{\mu_2}(\omega))+2\mu\kappa^2_{\mu_2}(\omega))\,dm\ge 0$. In particular, $\kappa_{\mu_2}^2$ is different from $0$ with respect to $m$. A symmetric procedure performed for $\mu<\mu_1$ (now defining $\kappa^1_\mu$ as the upper delimiter of the basin of attraction of the graph of $\alpha_\mu$, so that $\mu\mapsto\kappa^1_\mu$ is nonincreasing on $\mathcal I_1$) shows the existence of a $\tau_{\mu_1}$-equilibrium $\kappa_{\mu_1}^1\le 0$ with is different from 0 with respect to $m$. Finally, we assume for contradiction that $\mu_1=\mu_2$, observe that $\kappa^1_{\mu_1}\le 0\le \kappa^2_{\mu_2}$ define three $\tau_{\mu_1}$-equilibria which are different with respect to $m$, and conclude that $\int_\Omega f_x(\w,0)\,dm>0$ (see Subsection~\ref{subsec2:dconcavity}), which is not the case. This contradiction proves the assertion.

It remains to check that $\mathcal M_\mu^l=\mathcal M_0$ for $\mu>\mu_1$ and $\mathcal M_\mu^u=\mathcal M_0$ for $\mu<\mu_2$, which ensures that $\mathcal M_0$ is the unique $\tau_\mu$-minimal set for $\mu\in(\mu_1,\mu_2)$. We begin by assuming that $\mathcal M_{\mu_2}^u=\mathcal M_0$, which means that $\beta_{\mu_2}$ vanishes on the residual set of its continuity points. Proposition \ref{prop:7properties}(i) ensures that, if $\mu<\mu_2$, then $0\le\beta_\mu(\omega)\le\beta_{\mu_2}(\omega)=0$ at a common continuity point of both maps, and this ensures that $\mathcal M_\mu^u=\mathcal M_0$, as asserted. Now, let us work in the case $\mathcal M_{\mu_2}^u>\mathcal M_0$. The argument adapts to this situation that of the two last paragraphs of the proof of Theorem~\ref{th:clasification}, as we sketch in what follows. First, for $\mu>\mu_2$, we consider the $\tau_\mu$-equilibrium $\kappa_\mu^2$ of the previous paragraph, which satisfies $\kappa_\mu^2(\omega)>0$ $m$-a.e. if $\int_\Omega f_x(\omega,0)\, dm<0$ for an $m\in\mathfrak{M}_\mathrm{erg}(\Omega,\sigma)$. Second, we combine this property with the nonhyperbolicity of $\mathcal M^u_{\mu_2}$ to deduce that there are points in the graph of $\kappa^2_{\mu_2}=\lim_{\mu\downarrow\mu_2}\kappa^2_\mu$ which are between the delimiter equilibria of $\mathcal M_{\mu_2}^u$. And third, we deduce from this fact that $\beta_\mu(\omega)=0$ on its residual set of continuity points if $\mu<\mu_2$, which proves that $\mathcal M_\mu^u=\mathcal M_0$ also in this case. The argument is analogous for $\mathcal M_\mu^l$ and $\mu>\mu_1$, and the proof is complete.\end{proof}
Note that the model analyzed in Proposition~\ref{prop:4Bpextended} fits in the situation of Theorem~\ref{th:xsquareclasification}(iii), and that in that case we can determine the values of $\mu_1$ and $\mu_2$. Autonomous cases $x'=f(x)+\lambda x^2$ fitting the possibilities described in the previous theorem are very easy to find, since they just depend on the sign of $f'(0)$. For example, $x'=-x^3+x+\lambda x^2$ for (i), $x'=-x^3-x+\lambda x^2$ for (ii), and $x'=-x^3+\lambda x^2$ for (iii).

We close this paper by using the information just obtained to go deeper in the analysis of the bifurcation possibilities for our initial problem \eqref{eq:65}, i.e., $x'=f(\omega{\cdot}t,x)+\lambda x$. More precisely, Corollary \ref{coro:6finale} shows that for any $\lambda_0\le-\inf\text{\rm sp}(f_x(\cdot,0))$ there exists a suitable $\mu_0\in\mathbb R$ such that $\lambda_0$ is the lower bifurcation point of the modified family $x'=f(\omega{\cdot}t,x)+\mu_0x^2+\lambda x$, and that the three different possibilities of Theorem \ref{th:clasification} correspond to $\lambda_0<-\sup\text{\rm sp}(f_x(\cdot,0))$, $\lambda_0=-\inf\text{\rm sp}(f_x(\cdot,0))$ and $-\sup\text{\rm sp}(f_x(\cdot,0))\le \lambda_0<-\inf\text{\rm sp}(f_x(\cdot,0))$. To this end, we consider the two-parametric bifurcation problem of minimal sets
\begin{equation}\label{eq:7biparametric}
x'=f(\omega{\cdot}t,x)+\mu x^2+\lambda x\,,
\end{equation}
where $f\in C^{0,2}(\Omega\times\mathbb{R},\mathbb{R})$ is $\mathrm{(Co)}_2$ and $\mathrm{(SDC)}_*$ and $f(\omega,0)=0$ for all $\omega\in\Omega$. Let $\mathcal{A}_{\lambda,\mu}$ be the global attractor of the local skewproduct flow $\tau_{\lambda,\mu}$ induced on $\Omega\times\mathbb{R}$ by $\eqref{eq:7biparametric}_{\lambda,\mu}$, with delimiter equilibria $\alpha_{\lambda,\mu}(\omega)=\inf(\mathcal{A}_{\lambda,\mu})_\omega$ and $\beta_{\lambda,\mu}(\omega)=\sup(\mathcal{A}_{\lambda,\mu})_\omega$, and let $\mathrm{sp}(f_x(\cdot,0))=[-\lambda_+,-\lambda_-]$. We define
\begin{equation}\label{eq:7aplications}
\begin{split}
&\hat\mu\colon(-\infty,\lambda_+]\rightarrow\mathbb{R}\,,\;\lambda\mapsto\inf\{\mu\in\mathbb{R}\colon\,\text{the graph of }\beta_{\lambda,\xi}\text{ defines a}
\\&\qquad\qquad\qquad\qquad\text{hyperbolic minimal set }\mathcal{M}_{\lambda,\xi}^u\neq\mathcal{M}_0\text{ for all }\xi>\mu\}\,,\\
&\hat\lambda\colon\mathbb{R}\rightarrow(-\infty,\lambda_+]\,,\;\mu\mapsto\inf\{\lambda\in\mathbb{R}\colon\,\text{the graph of }\beta_{\xi,\mu}\text{ defines a}
\\&\qquad\qquad\qquad\qquad\text{hyperbolic minimal set }\mathcal{M}_{\xi,\mu}^u\neq\mathcal{M}_0\text{ for all }\xi>\lambda\}\,.
\end{split}
\end{equation}
Notice that, for any fixed $\lambda_1\le\lambda_+$, $\mathrm{sp}(f_x(\cdot,0)+\lambda_1)=[\lambda_1-\lambda_+,\lambda_1-\lambda_-]\not\subset(0,\infty)$, which ensures that the bifurcation diagram for the $\mu$-parametric family \eqref{eq:7biparametric}$_{\lambda_1}$ is not that of Theorem~\ref{th:xsquareclasification}(i) and hence that $\hat\mu(\lambda_1)$ is well defined. In addition, Theorem~\ref{th:clasification} ensures that $\hat\lambda$ is well defined and satisfies $\hat\lambda(\mu_1)\le\lambda_+$ in the three bifurcation cases which may arise for
the $\lambda$-parametric family \eqref{eq:7biparametric}$_{\mu_1}$, for any $\mu_1\in\mathbb R$.
\begin{proposition}\label{prop:6aplications} Let $\hat\lambda$ and $\hat\mu$ be the maps defined in \eqref{eq:7aplications}. Then,
\begin{enumerate}[label=\rm{(\roman*)}]
\item $\hat\lambda\circ\hat\mu=\mathrm{Id}_{(-\infty,\lambda_+]}$, and consequently, $\hat \lambda$ is onto.
\item $\hat\lambda$ is nonincreasing, and consequently, $\hat\lambda$ is continuous.
\end{enumerate}
\end{proposition}
\begin{proof} (i) We fix $\lambda_0\in(-\infty,\lambda_+]$ and call $\mu_0=\hat\mu(\lambda_0)$ and $\overline\lambda_0=\hat\lambda(\mu_0)$. The goal is to check that $\overline\lambda_0=\lambda_0$. The definition of $\mu_0$ ensures that the graph of $\beta_{\lambda_0,\mu_0}$ does not define a hyperbolic $\tau_{\lambda_0,\mu_0}$-minimal set distinct from $\mathcal{M}_0$, so $\lambda_0\le\overline\lambda_0$. By contradiction, we assume that $\lambda_0<\overline\lambda_0$, and fix $\lambda\in(\lambda_0,\overline\lambda_0)$. Observe that $\mu_0$ is the upper bifurcation point of the diagram described by Theorem~\ref{th:xsquareclasification}(ii) or (iii) for the $\mu$-family of \eqref{eq:7biparametric}$_{\lambda_0}$, and that, for any $\epsilon>0$, the upper delimiter $\beta_{\lambda_0,\mu_0+\epsilon}$ of $\mathcal{A}_{\lambda_0,\mu_0+\epsilon}$ is continuous and strictly positive. Then, if $\epsilon\in(0,1)$ satisfies $\epsilon<(\lambda-\lambda_0)/\sup_{\omega\in\Omega} \beta_{\lambda_0,\mu_0+1}(\omega)$, we have $\lambda_0-\lambda+\epsilon\beta_{\lambda_0,\mu_0+\epsilon}(\omega)<0$ for all $\omega\in\Omega$, and hence
\begin{equation*}
\begin{split}
\beta_{\lambda_0,\mu_0+\epsilon}'(\omega)&=f(\omega,\beta_{\lambda_0,\mu_0+\epsilon}(\omega))+\lambda_0\beta_{\lambda_0,\mu_0+\epsilon}(\omega)+(\mu_0+\epsilon)\beta_{\lambda_0,\mu_0+\epsilon}^2(\omega)\\
&= f(\omega,\beta_{\lambda_0,\mu_0+\epsilon}(\omega))+\lambda\beta_{\lambda_0,\mu_0+\epsilon}(\omega)+\mu_0\beta_{\lambda_0,\mu_0+\epsilon}^2(\omega)\\
&\qquad\qquad\qquad\qquad+\beta_{\lambda_0,\mu_0+\epsilon}(\omega)(\lambda_0-\lambda+\epsilon\beta_{\lambda_0,\mu_0+\epsilon}(\omega))\\
&< f(\omega,\beta_{\lambda_0,\mu_0+\epsilon}(\omega))+\lambda\beta_{\lambda_0,\mu_0+\epsilon}(\omega)+\mu_0\beta_{\lambda_0,\mu_0+\epsilon}^2(\omega)\,.
\end{split}
\end{equation*}
That is, $\beta_{\lambda_0,\mu_0+\epsilon}$ is a global strict lower solution for \eqref{eq:7biparametric}$_{\lambda,\mu_0}$. This ensures that $0<\beta_{\lambda_0,\mu_0+\epsilon}<\beta_{\lambda,\mu_0}$ (see Subsection~\ref{subsec:2compactcoerciveattractor}), which in particular implies the existence of a strictly positive $\tau_{\lambda,\mu_0}$-minimal set. But this is not possible: Theorem~\ref{th:clasification} and the definition of $\hat\lambda$ show that there is no $\tau_{\lambda,\mu_0}$-minimal set above $\mathcal{M}_0$ if $\lambda<\overline\lambda_0$. This contradiction shows that $\overline\lambda_0=\lambda_0$.

(ii) By contradiction, we assume that there exist $\mu_1<\mu_2$ such that $\lambda_1=\hat\lambda(\mu_1)<\hat\lambda(\mu_2)=\lambda_2$. We take $\lambda\in(\lambda_1,\lambda_2)$. As $\lambda>\lambda_1$, Theorem~\ref{th:clasification} and the definition of $\lambda_1$ ensure that the upper delimiter $\beta_{\lambda,\mu_1}$ of $\mathcal{A}_{\lambda,\mu_1}$ is continuous and strictly positive. Notice that
\begin{equation*}
\begin{split}
\beta'_{\lambda,\mu_1}(\omega)&=f(\omega,\beta_{\lambda,\mu_1}(\omega))+\lambda\beta_{\lambda,\mu_1}(\omega)+\mu_2\beta^2_{\lambda,\mu_1}(\omega)+(\mu_1-\mu_2)\beta^2_{\lambda,\mu_1}(\omega)\\
&<f(\omega,\beta_{\lambda,\mu_1}(\omega))+\lambda\beta_{\lambda,\mu_1}(\omega)+\mu_2\beta^2_{\lambda,\mu_1}(\omega)\,,
\end{split}
\end{equation*}
so $\beta_{\lambda,\mu_1}$ is a global strict lower solution for \eqref{eq:7biparametric}$_{\lambda,\mu_2}$. Therefore, $0<\beta_{\lambda,\mu_1}<\beta_{\lambda,\mu_2}$ (see Subsection~\ref{subsec:2compactcoerciveattractor}), which in particular implies the existence of a strictly positive $\tau_{\lambda,\mu_2}$-minimal set. But the definition of $\lambda_2$ ensures that $\mathcal{M}_0$ is the unique $\tau_{\lambda,\mu_2}$-minimal set, since $\lambda<\lambda_2$. Hence, $\hat\lambda$ is nonincreasing. Finally, a nondecreasing and onto function defined from an interval to an interval is always continuous.
\end{proof}
\begin{corollary} \label{coro:6finale}Let $\mathrm{sp}(f_x(\cdot,0))=[-\lambda_+,-\lambda_-]$ and $\lambda_0\le\lambda_+$. The $\lambda$-parametric family
\begin{equation}\label{eq:6endexplanation}
 x'=f(\omega{\cdot}t,x)+\hat\mu(\lambda_0) x^2+\lambda x
\end{equation}
exhibits
\begin{itemize}
\item[-] the local saddle-node and transcritical bifurcations of Theorem~{\rm \ref{th:clasification}(i)} if $\lambda_0<\lambda_-$, with $\lambda_0$ as local saddle-node bifurcation point and $\alpha_\lambda$ colliding with $0$ as $\lambda\downarrow\lambda_+$;
\item[-] the classical pitchfork bifurcation of Theorem~$\mathrm{\ref{th:clasification}(ii)}$ if $\lambda_0=\lambda_+$, with unique bifurcation point $\lambda_0$;
\item[-] and the generalized pitchfork bifurcation of Theorem~$\mathrm{\ref{th:clasification}(iii)}$ if $\lambda_0\in[\lambda_-,\lambda_+)$, with $\lambda_0$ as lower bifurcation point and $\alpha_\lambda$ colliding with $0$ as $\lambda\downarrow\lambda_+$.
\end{itemize}
\end{corollary}
\begin{proof} Proposition~\ref{prop:6aplications}(i) ensures that $\hat\lambda(\hat\mu(\lambda_0))=\lambda_0$. That is, $\lambda_0=\inf\{\lambda\in\mathbb{R}\colon$ the graph of $\beta_{\xi,\hat\mu(\lambda_0)}$ defines a hyperbolic minimal set $\mathcal{M}_{\xi,\hat\mu(\lambda_0)}^u\neq\Omega\times\{0\}$ for all $\xi>\lambda\}$. The conclusions follow from the description of the three cases made in Theorem~\ref{th:clasification} applied to \eqref{eq:6endexplanation}.
\end{proof}
By taking the lower delimiters of the global attractors instead of the upper ones in \eqref{eq:7aplications}, we get a result analogous to Corollary \ref{coro:6finale}, with $\beta_\lambda$ colliding with $0$ at the upper bifurcations points.

\end{document}